\magnification=\magstep1
\input amstex
\documentstyle{amsppt}
\TagsOnRight
\hsize=5.1in                                                  
\vsize=7.8in

\define\fp{\Bbb F_p}

\topmatter

\title Geometry of growth:\\
approximation theorems for $L^2$
invariants
\endtitle
\rightheadtext{Michael Farber}
\leftheadtext{Geometry of growth}
\author  Michael Farber \endauthor
\address
School of Mathematical Sciences,
Tel-Aviv University,
Ramat-Aviv 69978, Israel
\endaddress
\email farber\@math.tau.ac.il
\endemail
\thanks{The research was partially supported by a grant from the
Israel Academy of Sciences and Humanities}
\endthanks
\abstract{In this paper we study the problem of approximation of the
$L^2$-topological invariants by
their finite dimensional analogues. We obtain generalizations of the 
theorem of L\"uck \cite{L},
dealing with towers of finitely sheeted normal coverings. We
prove approximation theorems,
establishing relations between the homological invariants,
corresponding to infinite dimensional representations and sequences of
finite dimensional representations, assuming that their normalized characters converge.
Also, we find an approximation theorem for residually finite $p$-groups ($p$ is a prime),
where we use the homology with coefficients in a finite field $\fp$.}

We view sequences of finite dimensional flat bundles
of growing dimension as examples of growth processes. We study a
von Neumann category with a Dixmier type trace, which allows to describe the asymptotic
invariants of growth processes. We introduce a new invariant of torsion objects, the torsion
dimension. 
We show that the torsion dimension appears in general as an
additional correcting term in the approximation theorems; 
it vanishes under some arithmeticity assumptions.
We also show that the torsion dimension allows to establish non-triviality
of the Grothendieck group of torsion objects.

\endabstract
\endtopmatter
\footnote""{1991 Mathematics Subject Classification: Primary 53C23, Secondary
55N25}


\define\V{{\Cal V}}  
\define\C{{\Bbb C}}
\define\R{{\Bbb R}}  
\define\RR{{\Cal R}}       
\define\Z{{\Bbb Z}}  
\define\Q{{\Bbb Q}} 
  
\define\T{{\Cal T}}  
\define\X{{\Cal X}}   
\define\Y{{\Cal Y}}  
\define\M{{\Cal M}}    

\define\NN{\Cal N}
\define\N{{\Bbb N}}
 
\define\Hom{\operatorname{Hom}}

\define\Tr{\operatorname{Tr}}

\define\End{\operatorname{End}}   

\define\GL{\operatorname{GL}} 
        
\define\A{{\Cal A}}

\define\F{{\Cal F}}      
 
\define\tr{\operatorname{tr}} 
\define\im{\operatorname{im}} 
   

\define\id{\operatorname{id}}    

\define\ns{\frak {ns}}     
\define\ob{\operatorname{ob}}   
\define\cl{\operatorname{cl}} 
  
\define\rank{\operatorname{rank}}  
\redefine\H{\Cal H}

\define\E{\Cal E}
\define\U{{\Cal U}}
\def\<{\langle}
\def\>{\rangle}

\define\pd#1#2{\dfrac{\partial#1}{\partial#2}}
   
\define\cc{\Cal C}
  
\define\eca{\Cal E({\Cal C}_{\Cal A})}

\define\ca{\Cal C_{{\Cal A}}}

\define\tca{\T(\ca)}
\redefine\L{\Cal L}

\define\cm{{\Cal C}(\mu)}
\define\tdim {{\frak {tordim}}}
\define\Lim{\operatorname {Lim}}
\define\tm{\tr_{\omega}}
\define\dd{{\frak {projdim_\omega}}}
\define\lo{\Lim_{\omega}}

\define\lam{\lambda}
\define\Mat{\operatorname {Mat}}
\define\m{\frak m}
\define\ld{\operatorname{log\ det^\prime}} 
\define\W{\Cal W}
\define\p{\frak p}
\define\oo{\frak o}

\documentstyle{amsppt}   
\nopagenumbers

\heading{\bf \S 0. Introduction}\endheading

\subheading{0.1}
In 1994 Wolfgang L\"uck \cite{L} proved a beautiful theorem stating that 
von Neumann Betti numbers of the universal covering of a finite polyhedron
can be found as the limits of the
normalized Betti numbers of finitely sheeted normal coverings. Before L\"uck it was only
known that there is an inequality (called Kazhdan's inequality \cite{Ka}, cf. also
Gromov \cite{Gr}, pages 13 and 153).

One of the goals of the present paper is to generalize the L\"uck's theorem in two directions.
First, instead of finitely sheeted normal coverings we consider flat vector bundles
of finite dimension. Secondly, instead of $L^2$-Betti numbers we study the von Neumann dimensions
of the homology of infinite dimensional flat bundles determined by unitary representations
in a von Neumann category with a trace.

The other main purpose of this paper is to investigate the situations, when the statement
of the L\"uck's theorem in its original
form is incorrect. We show that the correcting additional term has a very interesting
meaning ({\it the torsion dimension});
it can be understood in the
framework of the formalism of extended cohomology and von Neumann categories.
As examples in the paper show,
vanishing of this correcting term happens in fact rarely, under very special arithmetic
assumptions.

\subheading{0.2} In order to illustrate our results,
we formulate here three approximation theorems, dealing with the towers of
coverings and the $L^2$-Betti numbers, which are corollaries of the main
Theorem 9.2 below. The first (Theorem 0.3) generalizes the L\"uck's theorem by admitting
towers of non-normal finitely sheeted coverings.
The second (Theorem 0.4) generalizes Theorem 0.3 by allowing twisted coefficients;
here we impose some important restrictions coming from the algebraic number theory.

\proclaim{0.3. Theorem} Let $\pi$ be an infinite discrete group and let
$\pi\supset \Gamma_1\supset\Gamma_2\supset\dots$ be a sequence of subgroups of finite index.
For any $k = 1, 2, \dots$ denote by $n_k$ the total number of the
subgroups of $\pi$, conjugate to $\Gamma_k$; given $g\in \pi$, we denote by $n_k(g)$ the
number of subgroups conjugate to $\Gamma_k$, containing $g$. Suppose, that for $g\ne 1$ holds
$$
\lim_{k\to\infty} \frac{n_k(g)}{n_k} = 0.\tag0-1
$$
For any finite polyhedron $X$ with $\pi_1(X)=\pi$, consider $[\pi:\Gamma_k]$-sheeted coverings
$\tilde X_k \to X$ corresponding to the subgroups $\Gamma_k \subset \pi$,
where $k =1, 2, \dots$. Then the sequence of the normalized Betti numbers
$$
\lim_{k\to \infty} \frac{\dim H_i(\tilde X_k)}{[\pi:\Gamma_k]} = b^{(2)}_i(X)\tag0-2
$$
converges to the $L^2$-Betti number $b^{(2)}_i(X)$.
\endproclaim

Proof of Theorem 0.3 will be given in \S 9.

Note that the condition (0-1) of Theorem 0.3 implies that $\pi$ is residually finite:
for fixed $k$ denote by
$P_k\subset \pi$ the intersection of all the subgroups of $\pi$, conjugate to $\Gamma_k$;
then because of (0-1) we have $\cap P_k = \{1\}$.

Assuming that all subgroups $\Gamma_k\subset \pi$ are normal, Theorem 0.3 reduces to the theorem
of L\"uck.

Here is another generalization of the theorem of L\"uck:

\proclaim{0.4. Theorem} Let $\pi$ be an infinite discrete group and let
$\pi\supset \Gamma_1\supset\Gamma_2\supset\dots$ be a sequence of normal
subgroups of finite index such that the intersection $\cap \Gamma_k = \{1\}$ is trivial.
Let
$\rho_k: \Gamma_k \to \Mat(m_k\times m_k, \oo_k)$ be a sequence of unitary representation,
where $k = 1, 2, \dots$; here $\oo_k$ denotes the ring of algebraic
integers in a number field $\F_k \subset \C$.
We assume that each $\F_k$ comes imbedded into $\C$ such that it is
invariant under the complex conjugation and we consider the induced involution on $\F_k$ and
on $\oo_k$. Let $h_k$ denote the degree of $\F_k$ over $\Q$; we assume that 
there exists a common upper bound
$$
h_k \le h.\tag0-3
$$
For any finite polyhedron $X$ with $\pi_1(X) = \pi$, consider the normal covering
$\tilde X_k\to X$, where $k=1, 2, \dots$,
corresponding to the subgroup $\Gamma_k$ and denote by $V^k$ the flat
bundle over $\tilde X_k$, determined by the representation $\rho_k$.
Then the sequence
$$
\lim_{k\to\infty} \ \frac{\dim H_i(\tilde X_k, V^k)}{\dim V^k\cdot [\pi:\Gamma_k]} \
= \ b^{(2)}_i(X)\tag0-4
$$
converges to the $L^2$-Betti number $b^{(2)}_i(X)$ of $X$.
\endproclaim

The proof of Theorem 0.4 will be given in 9.3; it follows from a more general Theorem 9.2,
dealing with sequences of flat bundles, satisfying some arithmeticity conditions, and
such that their normalized
characters converge to the character of a unitary representation in a von Neumann category.
In particular, the condition (0-3) is one of the main properties of {\it the arithmetic
approximation}, cf. 9.1. We show in \S 10 that Theorem 9.2 is false if the degrees $h_k$
are unbounded. The properties of arithmetic approximation include also an important condition
on the sequence of Galois groups acting on the characters of the approximating
sequence of representations; we show in \S 10 that the theorem becames false, if this condition
is violated.

It is interesting to emphasize that under the conditions of arithmetic approximation
the dimensions of the
flat bundles, approximating a von Neumann flat bundle, have to tend to infinity.

Here is another corollary of Theorem 9.2, which we prove in \S 9:

\proclaim{0.5. Theorem} Let $X$ be a finite polyhedron, and let $\rho: \pi_1(X) \to
\Mat(m\times m, \oo_{\F})$ be a unitary representation, where $\F\subset \C$ denotes a
cyclotomic field and $\oo_{\F} \subset \F$ denotes its ring of algebraic integers.
Suppose that $\rho$ is injective and its image has trivial intersection with the
the center of
the matrix algebra $\Mat(m\times m, \oo_{\F})$. Let $\E\to X$ denote the flat vector
bundle of rank $m$ determined by the representation $\rho$. Then 
$$
\lim_{k\to\infty} \ \frac{\dim H_i(X, \E^{\otimes k})}{(\dim \E)^k} \ = \ b^{(2)}_i(X).\tag0-5
$$
\endproclaim

A more general statements of this type can be found in \S 9, cf. Theorem 9.6.
 
\subheading{0.6} In this paper we use the language of
von Neumann categories, which provides a natural environment for developing the
$L^2$-homology theory, cf. \cite{F}. We review this material briefly in \S 1. Traces on
von Neumann categories play an important role; the traces allow to assign dimensions to
objects of the von Neumann category, which generalize the von Neumann dimension.

Given a polyhedron $X$ and a representation of the fundamental group $\pi$ of $X$
on an object of a von Neumann category with a trace, they determine {\it a character} on
$\pi$. It is a
class function $\chi: \pi \to \C$, which satisfies certain positivity condition, cf. \S 7.
We show here that knowing this character as the only information about the representation
allows to find the von Neumann Betti numbers and the spectral density function of the extended
$L^2$-homology. Conversely, we show that one may construct von Neumann categories with traces
starting from class functions on the fundamental group $\pi$.

The problem of describing the behavior of the $L^2$-invariants under deformations
of the von Neumann representation, seems to be of central importance. For example, one
wants to approximate von Neumann representations by finite dimensional
ones (as in the L\"uck's theorem).
Since the character of a von Neumann representation determines completely the
most important $L^2$-homological invariants, we study situations, when we have a
sequence of finite dimensional representations with the property that their normalized
characters converge (pointwise, i.e. as functions on the group $\pi$)
to the character of the given infinite dimensional representation. Our aim is to
find the homological (spectral) invariants corresponding to the infinite
dimensional representation in terms
of the approximating finite dimensional family; this seems to be
a natural generalization of the situation studied by L\"uck \cite{L}.

\subheading{0.7} It turns out that any approximating sequence of finite dimensional
representations
can also be treated as a single representation is a finite von Neumann category.
Moreover, this von Neumann category admits {\it a Dixmier type} (i.e. not normal)
{\it trace}; the construction of this trace
uses universal summation machines of von Neumann \cite{vN}. Note that Dixmier type traces
play a very important role in the noncommutative geometry of A. Connes \cite{C}.
We show in \S 2, that not normal
traces allow to define a dimension type function for the
torsion objects of the extended category. We call this function {\it the torsion dimension}.
Its main property is that it determines a non-trivial homomorphism on the Grothendieck group
of the torsion subcategory.

This von Neumann category
allows to study {\it the growth processes} - families of finite dimensional chain complexes.
A sequence of flat bundles over a finite polyhedron (more precisely, the corresponding
sequence of the chain complexes) is an instance of a growth process. Any growth process
defines its asymptotic invariants: the projective dimension, the torsion dimension, and the
spectral density function. As another geometrically interesting example of growth processes
we may mention the sequence of choppings (exhaustion) of a non-compact Riemannian manifold.
  
\subheading{0.8} In the most general approximation theorems established in \S 8
(cf. Theorems 8.2, 8.3), we find
that the torsion dimension of the extended homology appears as the additional correcting
term. In many cases one may expect the torsion dimension to be
independent of the choice of the summation machine $\omega$, which is a part of the Dixmier
type trace. We show that such independence happens in the analytic situation (Theorem 8.4). 
We also analyse examples showing that sometimes one may realize a sequence of approximating
Betti numbers by an arbitrary sequence consisting of 0's and 2's, cf. 6.3.

However, if we want
to guarantee vanishing of the torsion dimension in the general approximation theorem 8.2,
we have to impose some assumptions from algebraic number theory. The idea of {\it integrality}
is also very important in the original L\"uck's theorem. We develop this idea further, by
allowing representations over the algebraic integers of algebraic number fields;
this adds flexibility and makes possible many interesting applications.

\subheading{0.9} Finally, I want to mention an approximation theorem of a different type;
it is Theorem 11.1. Here we
assume that the fundamental group admits a chain of normal subgroups with index being a power
of $p$, where $p$ is a fixed prime number. We show that the Betti numbers over the finite field
$\fp$ behave {\it in a monotone fashion}; this produces an inequality
$$b_i^{(2)}(X) \le \dim_{\fp} H_i(X,\fp)\tag0-6$$
between the $L^2$-Betti numbers and the usual $\fp$-homology, cf. Corollary 11.2.

\subheading{0.10} This paper was written while the author was visiting IHES in
Bures-sur-Yvette (France); I am very thankful to the IHES for hospitality.

I am also grateful to D. Burghelea, A. Connes and M. Gromov for a number of
stimulating discussions.
\newpage

\heading{\bf \S 1. A review of extended $L^2$-homology, \\
von Neumann categories, and traces}\endheading

Intuitively, the extended homology provides a rigorous formalism to study 
a homology theory based on the (usual) infinite $L^2$-cycles together with the
"cycles" of the form
\vskip 4cm
\centerline{Fig. 1}
More precisely, we study geometry of non-compact manifolds or flat infinite dimensional
bundles over compact manifolds; the cycle on the figure 1 above represents
in fact a sequence of cycles $c_n$, where $n = 1, 2, \dots$
such that each $c_n$ is a boundary, but the {\it size} of a
minimal chain, spanned by
$c_n$ is much greater ({\it asymptotically}) than the volume of $c_n$. 

A precise definition of the extended $L^2$ homology uses a generalization
of the notion of Hilbert space - the functor of extended homology assigns to a manifold such
generalized Hilbert space. It turns out that the familiar category of Hilbert spaces is not
good enough; we complete it by adding {\it "torsion Hilbert spaces"},
such that the obtained category becomes an abelian category.
In order to obtain a good category and to include some interesting applications,
it is reasonable to study this construction
of abelian extension starting from a von Neumann category.

In this section we will give a brief review of the notion of von Neumann category,
the extended abelian categories, and traces, which will be used in the rest of the paper.
In full detail all this material is described in \cite{F}.

\subheading{1.1. von Neumann categories}Let $\A$ be an algebra over $\C$ having an involution 
which will be denoted 
by the star $\ast$.  {\it A Hilbert representation of $\A$ (or a Hilbert module)} 
is a Hilbert space $\H$ supplied
with a left action of $\A$ on $\H$ by bounded linear maps such that for any
$a\in \A$ holds 
$$\<ax,y\>=\<x,a^\ast y\>\tag1-1$$ 
for all $x,y\in \H$. {\it A morphism} between Hilbert representations 
$\phi: \H_1\to \H_2$ is a bounded linear map commuting with the action of 
the algebra $\A$. We obtain the additive category of all Hilbert representations of a given
$\ast$-algebra $\A$. 

Assume that $\cc_\A$ is an additive subcategory of the category of 
all Hilbert representations of $\A$. We say that $\ca$ is {\it a von Neumann category} if
the following properties are verified:
\roster 
\item"{(i)}" {\it The kernel of any morphism $\phi:\H_1\to \H_2$ in
$\cc_\A$  and the natural inclusion $\ker \phi\to \H_1$ belong to $\cc_\A$.
\item"{(ii)}" For any morphism $\phi:\H_1\to \H_2$ of $\ca$ the adjoint
operator $\phi^\ast:\H_2\to \H_1$ is also a morphism of $\ca$.}
\item"{(iii)}" {\it for any pair of representations $\H_1, \H_2\in\ob(\ca)$, 
the corresponding 
set of morphisms $\Hom_{\ca}(\H_1,\H_2)$ is a weakly closed subspace 
in the space of all bounded linear
operators between $\H_1$ and $\H_2$.}
\endroster

Note, that for any object $\H\in\ob(\ca)$ of a von Neumann category the set of endomorphisms
$\Hom_{\ca}(\H,\H)$ is a von Neumann algebra.

\subheading{1.2. Finite objects}
We will say that an object $\H\in\ob(\ca)$ of a von Neumann category 
{\it is finite} if any closed 
$\ca$-submodule $\H_1\subset \H$ which is isomorphic to $\H$ in 
$\ca$, coincides with $\H$.

This property is equivalent to the requirement that
the von Neumann algebra $\Hom_{\ca}(\H,\H)$ of endomorphisms of $\H$
is finite. Cf. \cite{Di}, part III, chapter 8, \S 1.

A von Neumann category $\ca$ is called {\it finite} if all its objects are
finite.

\subheading{1.3. Trace and dimension} Let $\ca$ be a 
von Neumann category. 

{\bf Definition.} {\it A trace on category $\ca$ is a function, denoted $\tr$, 
which assigns to each object $\H\in\ob(\ca)$ a finite, non-negative trace 
$$\tr_{\H}: \Hom_{\ca}(\H,\H)\to \C\tag1-2$$
on the von Neumann algebra $ \Hom_{\ca}(\H,\H)$; in other words $\tr_{\H}$ assumes
(finite) values in
$\C$, $\tr_{\H}(a)$ is non-negative on positive elements $a$ of $ \Hom_{\ca}(\H,\H)$,
and $\tr_{\H}$ is traceful, i.e. 
$\tr_{\H}(ab)=\tr_{\H}(ba)$, for $a, b \in \Hom_{\ca}(\H,\H)$. 
It is also assumed that for any pair of representations
$\H_1$ and $\H_2$ the corresponding traces $\tr_{\H_1}$, $\tr_{\H_2}$ and 
$\tr_{\H_1\oplus\H_2}$ 
are related as follows: if 
$f\in\Hom_{\ca}(\H_1\oplus \H_2, \H_1\oplus \H_2)$ is given by a $2\times 2$
matrix$(f_{ij})$, where $f_{ij}:\H_i\to\H_j$, $i,j=1,2$, then}
$$\tr_{\H_1\oplus\H_2}(f)\ =\ \tr_{\H_1}(f_{11}) + \tr_{\H_2}(f_{22}).\tag1-3$$

For the notion of positive elements of the von Neumann algebra $\Hom_{\ca}(\H,\H)$
we refer to \cite{T}, page 24.

We will say that a trace $\tr$ on a von Neumann category is {\it normal} 
iff for each non-zero $\H\in\ob(\ca)$ the trace $\tr_\H$ on the von Neumann algebra
$\Hom_{\ca}(\H,\H)$ is normal. Recall that this means that
$$\sup_i \{\tr_\H(a_i)\} = \tr_\H(\sup_i\{ a_i\})\tag1-4$$
for any bounded increasing net
$a_i\in \Hom_{\ca}(\H,\H)$ consisting of positive operators; cf. \cite{T}, page 309. 

Given a trace $\tr$ on a category $\ca$, one can define the following {\it dimension function}:
$$\dim \H = \dim_{\tr} \H = \tr_{\H}(\id_{\H})\tag1-5$$
The real number $\dim_{\tr} \H$ is called {\it the von Neumann dimension
(or the projective dimension) of $\H$ with respect to the trace $\tr$}.

\subheading{1.4. The abelian extension} Given a von Neumann category 
$\ca$, there exists a bigger category $\E(\cc_\A)$, which is abelian and which
contains $\ca$ as a full subcategory. The construction of $\eca$ was suggested in \cite{F1},
\cite{F} using ideas of P. Freyd \cite{Fr}.

{\it An object} of the category $\E(\cc_\A)$ is defined as a
morphism $(\alpha:A^\prime\to A)$ in the category $\cc_\A$. 
Given a pair of objects $\X=(\alpha: A^\prime \to A)$ and 
$\Y=(\beta:B^\prime\to B)$ of $\E(\cc_\A)$,
a {\it morphism} $\X\to\Y$ in the category $\E(\cc_\A)$ is an  
equivalence class of morphisms $f:A\to B$ of category $\cc_\A$
such that $f\circ\alpha=\beta\circ g$
for some morphism $g:A^\prime \to B^\prime$ in $\cc_\A$. 
Two morphisms $f:A\to B$ and $f^\prime:A\to B$ of $\cc_\A$ represent 
{\it identical morphisms $\X\to\Y$ of $\E(\cc_\A)$} iff 
$f-f^\prime = \beta\circ F$ for some morphism $F:A\to B^\prime$ of category
$\cc_\A$. This defines an equivalence relation.
The morphism $\X\to\Y$, represented by $f:A\to B$,
is denoted by
$$[f]:(\alpha:A^\prime\to A)\ \to\ (\beta:B^\prime\to B)\quad
\text{or by}\quad
[f]:\X\to\Y.\tag1-6$$
The {\it composition} of morphisms is defined as the composition of the 
corresponding morphisms $f$ in the category $\cc_\A$.

\subheading{1.5. Embedding of $\ca$ into $\eca$} Given an object $A\in\ob(\ca)$
one defines the following object $(0\to A)\in\ob(\eca)$ of the extended
category. Since any morphism $f:A\to B$ determines a morphism 
$[f]:(0\to A)\to (0\to B)$ in the extended category, we obtain a full embedding $\ca\to \eca$.

It is possible to characterize the objects of the extended category which are
isomorphic in $\eca$ to objects coming from $\ca$ in intrinsic terms. Namely,
{\it an object $\X\in\ob(\eca)$ is projective if and
only if it is isomorphic in $\eca$ to an object of the form $(0\to A)$,
where $A\in\ob(\ca)$}

\subheading{1.6. The torsion subcategory}An object $\X=(\alpha:A'\to A)$ of the extended
category $\eca$ is called {\it torsion} iff the image of $\alpha$ is dense in $A$.

We will denote by $\tca$ the full subcategory of $\eca$ generated by all
torsion objects. $\tca$ is called {\it the torsion subcategory of $\eca$}.
If $\ca$ is a finite von Neumann category, then the torsion subcategory $\tca$ is an abelian 
subcategory of $\eca$.

Given an arbitrary object
$\X=(\alpha:A'\to A)$ of $\eca$ one considers the following torsion
object $T(\X) =\ (\alpha: A'\to \cl(\im(\alpha)))$
which is called {\it the torsion part of $\X$.} There is an obvious monomorphism
$T(\X)\to\X$. The factor $P(\X) = \X/T(\X)$ is projective,
called {\it the projective part of} $\X$.
We have
$\X = T(\X) \oplus P(\X).$
Thus, {\it the isomorphism type of an object of the extended category $\eca$
is determined by the isomorphism types of its projective and torsion parts.}

\subheading{1.7. Novikov - Shubin invariants} Given a trace on a von Neumann category
$\ca$, one obtains the numerical invariant $\ns(\X)$ of torsion objects,
called the Novikov - Shubin
invariant. We refer to \cite{F}, section 3.9, where it is described. There exist also
other invariants of torsion objects,
independent of the Novikov - Shubin invariant, cf. \cite{F1}. 

In the next section we will define new numerical invariant of torsion objects, which is
sometimes more convenient.

\subheading{1.8. Extended homology}
The functor of extended homology is constructed as follows, cf. \cite{F}, \cite{F1}.
Suppose that $X$ is a finite polyhedron with fundamental group $\pi$. Let $\ca$ be a
von Neumann category, and let $\rho: \pi\to \Hom_{\ca}(\M,\M)$ be a representation, where
$\M\in\ob(\ca)$. Consider the chain complex $C_\ast(\tilde X)$ (the cellular chain complex
of the universal covering $\tilde X$). Then
$$\M \otimes_\pi C_\ast(\tilde X)$$
is a chain complex in category $\ca$. Thus, it lies in the abelian category $\eca$)
and its homology (calculated in $\eca$,
called {\it extended $L^2$ homology of $X$ with coefficients in $\M$}) is
denoted by $\H_\ast(X,\M)$. Being an object of $\eca$, it is a direct sum of its projective
and torsion parts. The projective part of the extended homology coincides with the
reduced $L^2$ homology, cf. \cite{A} (defined by dividing the space of infinite $L^2$
chains by the closure
of $L^2$ boundaries). The torsion part of the extended homology is responsible for the
"almost cycles" or "asymptotic cycles" as the one shown on Figure 1.

\heading{\bf \S 2. Torsion dimension}\endheading

In this section we define a new numerical invariant of torsion objects, which we call 
torsion dimension. It behaves in better way, than
the known invariants (such as the Novikov-Shubin invariant and the minimal number of generators,
introduced in \cite{F}). We will use the torsion dimension in the next section to study the
Grothendieck group of torsion objects. Also, we will use the torsion dimension in approximation
theorems for $L^2$ topological invariants, cf. Theorems 8.2 and 8.4, where it produces a
correcting additional term. 

Everywhere in this section $\ca$ will denote a finite von Neumann category.
We will assume that we have a fixed trace $\tr$ on $\ca$, cf. subsection 1.3.
We will {\it not} assume that the trace $\tr$ is normal, since in the case of normal traces
the torsion dimension is always zero. Also, the trace $\tr$ is not supposed to be faithful.

Not normal traces are usually called {\it Dixmier type traces},
cf. \cite{C}, since J. Dixmier \cite{D} was the first who constructed such traces. Dixmier
type traces play very important role in the non-commutative geometry of A. Connes \cite{C}.

\subheading{2.1} First we will show that any non-normal trace determines {\it a dimension
function} of torsion objects. We will see that it behaves sub-additively under extensions.

Let $\X=(\alpha: A'\to A)$ be a torsion object of the extended category $\eca$
and let $F(\lambda)$ be its spectral density
function with respect to the trace $\tr$, cf. \cite{F}, formula (3-12).

\subheading{Definition} We will define {\it the torsion dimension of $\X$
(with respect to the trace $\tr$)} as the following real number
$$
\tdim \X = \tdim_{\tr} \X  = \lim_{\lambda\to +0} F(\lambda).\tag2-1
$$
Note that $F(\lambda)$ is increasing and so the limit exists.

We will also define the {\it reduced spectral density function} by
$$\tilde F(\lambda) = F(\lambda) - \tdim \X.\tag2-2$$

Note, that if the trace $\tr$ is normal, then the torsion dimension $\tdim \X$ is always zero.

\proclaim{2.2. Proposition} The torsion dimension $\tdim \X$ depends only on the isomorphism
type of $\X$ as an object of the extended category. The reduced spectral density functions
corresponding to isomorphic
torsion objects $\X$ and $\Y$ are dilatationally equivalent. \endproclaim

\demo{Proof} The proof of Proposition 3.8 in \cite{F} does not use the assumption of
normality of the trace. It shows that if $\X$ and $\Y$ are isomorphic torsion objects of the
extended category then the corresponding spectral density functions are dilatationally
equivalent. This implies our statement. \qed
\enddemo

Now we will establish the following internal characterization of the torsion dimension.
Let us recall that any trace $\tr$ on von Neumann category $\ca$ determines a dimension
function on $\ca$, cf. 1.3.

\proclaim{2.3. Proposition} Given a torsion object $\X$ of the extended category $\eca$,
its torsion dimension $\tdim \X$, equals to the
infimum of the von Neumann dimensions $\dim P$ (with respect to the trace $\tr$) of projective
objects $P$ of $\ca$ such that there exists an epimorphism $P\to \X$. \endproclaim

\demo{Proof} Suppose that $\X = (\alpha: A' \to A)$ and $\alpha$ is injective.
Recall that the spectral density function
$F(\lambda)$ is defined as follows. We consider the positive square root $T$ of the equation
$T^2 = \alpha^\ast\alpha$ and the spectral decomposition $T = \int_0^\infty \lambda dE_\lambda$.
Then $F(\lambda)$ is the von Neumann dimension of the subspace $E_\lambda A'$.

Thus, for any $\lambda >0$ the spectral projection $E_\lambda$ determines a projective
object $P=E_\lambda A'$ which has von Neumann dimension $F(\lambda)$
and which maps epimorphically onto $\X$. Indeed, the torsion object
$\X_\lambda = (\alpha: E_\lambda A' \to \alpha(E_\lambda A'))$ is isomorphic to $\X$ and
$E_\lambda A'$ and $\alpha(E_\lambda A')$ are isomorphic (by Lemma 2.3 of \cite{F}.
Therefore $\dim\X \ge \inf P$.

On the other hand, if $P$ is projective and maps epimorphically onto $\X$ then $\X$ admits
a representation of the form $(\gamma: P' \to P)$ with some $P'$ and
$\gamma$ and thus we obtain (using Proposition 2.2) that
$\tdim \X \le \inf P$. \qed \enddemo

Now we will show that the torsion dimension is sub-additive for extensions.

\proclaim{2.4. Proposition} For any short exact sequence
$$0\to \X' \to \X \to \X'' \to 0,\tag2-3$$
consisting of torsion objects of the extended abelian category $\eca$, holds
$$\max \{ \tdim(\X'), \tdim(\X'')\}  \le \tdim (\X) \le \tdim (\X') + \tdim (\X'').\tag2-4$$
Moreover, if the sequence (2-3) splits, then
$$\tdim (\X) = \tdim (\X') + \tdim (\X'').\tag2-5$$
\endproclaim

\demo{Proof} We will use the internal characterization of the torsion
dimension given by Proposition 2.3. It is clear that if $P'$ can be
mapped epimorphically onto $\X'$ and $P''$ can be mapped epimorphically
onto $\X''$, then their direct sum $P' \oplus P''$ can be mapped
epimorphically onto $\X$. Thus we obtain the right side of inequality (2-4).

>From Proposition 2.3 clearly follows that $\tdim (\X) \ge \tdim (\X'')$. Suppose now that
$P \to \X$ is an epimorphism with $P$ being a projective object of $\ca$. Let $P' \to P$
be the kernel of the composite $P \to \X \to \X''$. Then we have an epimorphism $P' \to \X'$.
Observe that $P'$ is isomorphic to $P$ in $\ca$ by Lemma 2.3 of \cite{F}; therefore
$\dim P = \dim P'$. This proves that $\tdim (\X) \ge \tdim (\X')$.

The equality (2-5) obviously follows from the definitions. \qed
\enddemo

\heading{\bf \S 3. Grothendieck group of torsion objects}\endheading

Note that equality (2-5) represents a very important distinction between the properties of two
functions on isomorphism types of torsion objects - the
torsion dimension, which we introduced above in \S 2, and the well known Novikov - Shubin
invariant. Recall that the Novikov - Shubin invariant of a direct sum equals to the minimum
of the Novikov - Shubin invariants of the summands:
$$\ns(\X' \oplus \X'') = \min \{\ns(\X'), \ns(\X'')\}.\tag3-1$$
The advantage of (2-5) is that it implies that {\it the torsion dimension determines
a homomorphism with values in $\R$ from the Grothendieck group constructed out of
abelian category $\tca$ of torsion objects in} $\ca$. Thus existence of a non-normal trace
on $\ca$ implies non-triviality of the Grothendieck group. We will make all
this precise in the following subsection.

\subheading{3.1. Grothendieck group of torsion objects of $\ca$} We will denote by
$K(\tca)$ the Grothendieck group of the abelian category $\tca$, cf. \cite{K}, page 53.
Recall that $K(\tca)$ is an abelian group generated by the symbols $[\X]$, one for each
isomorphism type of torsion objects $\X$ in $\ca$, with the addition given by
$[\X] + [\Y] = [\X\oplus\Y]$. The torsion dimension gives a well-defined homomorphism
$$\tdim: \ K(\tca) \to \R\tag3-2$$
(by Propositions 2.2 and 2.4).

\proclaim{3.2. Theorem} If the given trace $\tr$ on the category $\ca$ is non-normal, then
homomorphism (3-2) is non-trivial and thus the Grothendieck group $K(\tca)$ is non-zero.
\endproclaim
\demo{Proof} If the trace is non-normal then we may find a sequence
$$\H = \H_1 \supset \H_2 \supset \H_3 \supset \dots,$$
where $\H$ is an object of $\ca$ and $\H_n$'s are its closed subobjects, such that
$\cap \H_n = 0$ and $\lim \dim (\H_n) = c >0$. Define the following projector valued
function
$E_\lam$ for $\lam \in [0,1]$, by setting $E_\lam$ = the projection onto $\H_n$, for
$(n+1)^{-1} <\lam \le n^{-1}$. Then we consider the morphism
$$\alpha: \H \to \H,\quad\text{where}\quad \alpha = \int_0^1 \lambda dE_\lam.$$
Then $\X = (\alpha: \H \to \H)$ is a torsion object, and clearly $\tdim(\X) = c >0$.\qed
\enddemo

\heading {\bf \S 4. An example of von Neumann category with Dixmier type trace}\endheading

Our purpose now is to describe the simplest example of a finite
von Neumann category with a Dixmier type trace. This category will be important for our
applications to the problem of approximation of $L^2$ invariants; we will see in subsection
4.8 that this category allows to describe geometry and topology of growth processes.

\subheading{4.1} Fix a sequence of non-negative real numbers
$\mu = (\mu^n)$, $\mu^n > 0$, $\mu^n\in\R$. We will call $\mu$ {\it the growth rate}.
Normally, we will have in our applications $\mu^n$ tending to 0, or being constant.
The von Neumann category, we are going to construct will depend on this choice; we will denote
it $\cm$.

{\it Objects} of the category $\cm$ are sequences $\V = (V^n)$ of finite dimensional Euclidean
spaces, where $n$ runs over non-negative integers, such that the growth rate of the
dimension of $V^n$ is bounded above by the given sequence $\mu$:
$$ \dim V^n = O((\mu^n)^{-1}).\tag4-1$$
In other words, we assume that the product $\mu^n \cdot\dim V^n$ is bounded.
Note that each $V^n$ {\it is Euclidean}, i.e. it is supplied with a scalar product.

Each object $\V$ of $\cm$ determines {\it a Hilbert space} $\H_\V$, where
$$\H_{\V} = \{v = (v^n), v^n \in V^n; \sum ||v^n||^2 < \infty\}.\tag4-2$$
Here the norm $||v^n||$ denotes the norm of the space $V^n$.

{\it A morphism } $f: \V\to \W$ in $\cm$, where $\V=(V^n)$ and $\W=(W^n)$, is a sequence
$f=(f^n)$, where $f^n: V^n\to W^n$ is a linear map, such that there exists a common
upper bound
$$||f^n|| \le M\tag4-3$$
($M$ is independent of $n$).
Any morphism $f: \V\to \W$ of $\cm$
clearly induces a bounded linear map of the corresponding Hilbert
spaces, which we denote $f: \H_\V \to \H_\W$. Now one checks easily, that all properties from
the definition of von Neumann category (cf. 1.1) are satisfied. The algebra $\A$
in this case is $\A = \C$.

The category $\cm$ is clearly {\it finite}, cf. 1.2.

\subheading{4.2} Now we will describe a Dixmier type trace on $\cm$.
First, we recall from \cite{D} and \cite{C}, page 305, that there exists a linear form
$\lo$ (invented by J. von Neumann \cite{vN})
on the space $\ell^\infty(\N)$ of bounded sequences of complex numbers, that satisfies
the following conditions:

$$
\aligned
&(\alpha)\qquad \lo (\alpha_n) \ge 0 \quad\text{if}\quad \alpha_n\ge 0,\\
&(\beta)\qquad \lo (\alpha_n) = \Lim \alpha_n\quad\text{if}\quad
\alpha_n\quad\text{is convergent},\\
&(\gamma)\qquad \lo (\alpha_1,\alpha_1,\alpha_2,\alpha_2,\alpha_3,\alpha_3,\dots ) =
\lo (\alpha_n).
\endaligned
$$
Note that the form $\lo$ is not unique; it depends on the choice of the "rule" $\omega$, which
is sometimes called {\it the summation machine}.

Now for any object $\V$ of $\cm$ and for any endomorphism $f: \V\to \V$ in $\cm$ define its
trace $\tr_\omega(f)$ by
$$\tr_\omega(f) = \lo (\mu^n \cdot \Tr(f^n)).\tag4-4$$
Here on the right hand side of (4-4) $\Tr(f^n)$ denotes the usual finite dimensional
trace of the
linear map $f^n : V^n \to V^n$. Note that because of condition (4-3) we have
$|\Tr(f^n)| \le M \dim V^n$ and using (4-1) we see that the sequence $\mu^n \Tr(f^n)$
is bounded, and therefore the definition is correct.

It is easy to see that (4-4) defines a trace on the category $\cm$ (in the sense of
\cite{F}, Definition 2.7)
which is non-negative and traceful. We will see later that $\tm$ is not normal and not faithful.
Note that the constructed trace $\tm$ is not unique -- it depends on the choice of the
functional $\lo$ (i.e. on the "rule" $\omega$).

According to philosophy of A. Connes \cite{C},
in problems, having geometric origin, the answer will be often independent of $\omega$;
such problems A. Connes calls {\it measurable}. Cf. for example Proposition 5 in \cite{C},
chapter IV, section 2.$\beta$ concerning the Wodzicki residue.

We will also see examples of measurable problems (Theorems 8.4, 9.2 and 11.1) and
not measurable problems (example 6.3) later in this paper.

\subheading{4.3. The projective dimension} We know that any trace on a von Neumann
category determines a dimension function, cf. 1.3 above. The trace $\tr_\omega$ on
$\cm$ defined by (4-4) determines the following dimension function
$$\dd \V = \lo (\mu^n \cdot\dim V^n),\tag4-5$$
which we will call {\it the projective dimension} of $V$.

Note that the projective dimension $\dd \V$ depends only on the asymptotic behavior of the
numbers $\dim V^n$ for large $n$ and does not depend on any finite number of $\dim V^n$.
In particular, the
projective dimension $\dd \V$ vanishes if $V^n$ is non-zero only for finitely many $n$.
This shows
that the projective dimension $\dd$ (or, more precisely, the trace (4-4))
is not faithful -- nontrivial object may have trivial dimension.

Also, given an object $\V$ of $\cm$ with $\dd \V \ne 0$, consider the following sequence
$\V(m),\quad m=1,2,3,\dots$ of truncated objects of $\cm$,
where $\V(m)^n$ equals to $V^n$ for $n\le m$ and $V(m)^n = 0$ for $n > m$. We see that
$\V(m) \subset \V$ and
$$\sup_m \V(m) = \V.$$
However $\dd \V(m) = 0 \ne \dd \V$ for all $m$.
Therefore {\it the trace (4-4) is not normal}.

\subheading{4.4. The torsion dimension} Any torsion object of the extended
abelian category $\E(\cm)$, constructed out of $\cm$, cf. 1.4 and also \cite{F}, \S 1,
is represented by a morphism of $\cm$
$\X = (\alpha: \V \to \V),$.
Recall that $\alpha = (\alpha^n)$, where $\alpha^n : V^n \to V^n$ is a linear map.
We want to translate the general definition of the spectral density function
cf. \cite{F}, subsection 3.7, to the present
situation. Given a positive $\lambda >0$, denote by $F^n(\lambda)$ the maximal dimension of a
linear subspace contained in the following cone
$$\{v\in V^n; (\alpha^n(v),\alpha^n(v)) \le \lambda^2 (v,v)\}.\tag4-6$$
It is also equal to the number of eigenvalues of
$(\alpha^n)^\ast\alpha^n: V^n \to V^n$ which are
less than $\lambda^2$.
Then the spectral density function of $\X$ is given by
$$F(\lambda) = \lo (\mu^n F^n(\lambda)).\tag4-7$$
The {\it torsion dimension} (defined in \S 1 above) of $\X$ is by definition
$$\tdim_\omega \X = \lim_{\lambda\to +0} F^n(\lambda).\tag4-8$$
Roughly, the torsion dimension in this situation can be characterized as
{\it the density near zero of eigenvalues of $(\alpha^n)^\ast\alpha^n$
with respect to the chosen scale $\mu = (\mu^n)$}. 

\subheading{4.5. An example} Fix an arbitrary sequence $a^n$ of positive real numbers with tends
to 0. Consider an arbitrary object $\V = (V^n)$ of $\cm$ and let the morphism
$\alpha: \V\to \V$ be given as follows: $\alpha^n: V^n \to V^n$ is multiplication by $a^n$.
Then we obtain a torsion object $\X = (\alpha: \V\to \V)$ of $\cm$. For this $\X$ we have
$$
F^n(\lambda) =
\cases
0, \quad\text{if}\quad a^n > \lambda,\\
\dim V^n, \quad\text{if}\quad  a^n \le \lambda.
\endcases
$$
Thus we obtain that the spectral density function of $\X$
$F(\lambda) = \lo \mu^n F^n(\lambda)$
{\it is constant and equals the Dixmier dimension of} $\V$,
$F(\lambda) = \dd (\V).$
Therefore the torsion dimension of $\X$ (defined in \S 2) equals to $\dd(\V)$.

This is example shows that the torsion dimension may assume arbitrary non-negative real numbers.

\subheading{4.6. Extended homology in $\cm$} Consider a chain complex in $\cm$ of 
length $m$. Any such chain complex $C$ is just a sequence $C = (C^n, d^n)$, where $n=1,2,\dots$
of finite-dimensional complexes
$$C^n = (0\to C^n_m @>{d^n}>> C^n_{m-1} @>{d^n}>>\dots C^n_0\to 0)\tag4-9$$
such that

{\it (1) each chain space $C^n_i$, where $i = 0,1,\dots, m$, has a fixed Euclidean structure;

(2) the dimension growth rate satisfies $\dim C^n_i = O((\mu^n)^{-1})$;

(3) the norm of the differentials $d^n$ has a common upper bound $||d^n||\le M$.}

Given such chain complex $C$, it determines the extended homology, having the projective and
torsion parts, and we want to understand the Dixmier dimension of the projective part
and also the torsion dimension of the torsion part. Let
$Z_i^n$ denote
the space of cycles $\ker[d^n: C_i^n \to C_{i-1}^n]$; then $Z_i = (Z_i^n)$ and
$C_i = (C_i^n)$ are objects of
$\cm$. We clearly have for the extended homology of $C$:
$$\H_i(C) = (d: C_{i+1} \to Z_i).\tag4-10$$
The projective part of $\H_i(C)$ is just $(H_i(C^n))$, i.e. it is given by the sequence
consisting of the usual homology of the complexes $C^n$. Therefore, the Dixmier dimension
of the projective
part of the extended homology is given by
$$\dd (P(\H_i(C))) = \lo (\mu^n\cdot \dim H_i(C^n)).\tag4-11$$
Let $B_i^n$ be the subspace of boundaries $B_i^n = \im[d^n: C_{i+1}^n \to C_i^n]$
and let $B_i = (B_i^n)\in \ob(\cm)$. Then the
torsion part of the extended homology is given by
$$ T(\H_i(C)) = (d: C_{i+1}/Z_{i+1} \to B_i)\tag4-12$$

We summarize now the above discussion as follows:
\proclaim{4.7. Proposition} Suppose that a chain complex (4-9) in $\cm$ is given.
For any pair of integers
$n$ and $i$ denote by $h^n_i$ the number of zero eigenvalues of the operator
$$(d^n)^\ast d^n : C^n_{i+1} \to C^n_{i+1}\tag4-13$$
(the "Half-Laplacian") and for any $\lambda >0$ denote by $G^n_i(\lambda)$ the number of
eigenvalues
of (4-13) lying in the interval $(0,\lambda^2)$. Then the Dixmier dimension of the projective
part of the extended $i$-dimensional homology $\H_i(C)$ equals
$\lo \mu^n h^n_i$ and the spectral density function of the torsion part of $\H_i(C)$
is $G_i(\lambda) = \lo \mu^n G^n_i(\lambda)$. In particular, the torsion dimension of the
torsion part of $\H_i(C)$ is
$$\tdim_\omega T(\H_i(C)) = \lim_{\lambda\to +0} (\lo \mu^n G^n_i(\lambda))\tag4-14$$
\endproclaim

\subheading{4.8. Asymptotic invariants of a growth process}
A typical geometric situation, when the above numerical invariants of chain complexes in $\cm$
(the projective and torsion dimension and the Novikov - Shubin invariants)
can be applied consists in the following.

Suppose that $K =(K^n)$ is {\it growth process}, i.e.
a sequence of finite simplicial complexes, such that for any
integer $i$ the number of $i$-dimensional simplices in $K^n$ is $O((\mu^n)^{-1})$.

As a concrete
example (which will be studied in detail later in this paper)
we may assume that the complexes $K = (K^n)$ form a tower of finitely sheeted coverings
over a fixed finite polyhedron.

Another source of examples of growth processes is the following.
Suppose that we have an infinite polyhedron and the
finite polyhedra $K^n$ with $K^n\subset K^{n+1}$ form its exhaustion.
\vskip 4cm
\centerline{Fig. 2}
Growth process of this type
was considered in a recent preprint \cite{DM} of J. Dodziuk and V. Mathai.

Another example of a growth process provides a sequence of smaller and smaller polyhedral
approximations to a given compact Riemannian manifold.

Let us return now to the general situation. Given a growth process $(K^n)$,
we obtain the chain complexes $C^n = C_\ast(K^n)$ corresponding to the given
simplicial structures on the complexes $K^n$.
We may introduce the euclidean structure on $C_\ast(K^n)$, such that the simplices of
$K^n$ form an
orthonormal base. (Note, that in fact there may be different geometrically interesting ways of
choosing the scalar product on the chain space $C_\ast(K^n)$.)

We also have to verify condition (3) in subsection 4.6. Note that this
condition will be automatically satisfied if the growth process $(K^n)$ has
{\it bounded geometry:

there is a constant $M$ (independent of $n$)
such that each $i$-dimensional simplex of $K^n$
is adjacent in $K^n$ to at most $M$ simplices of dimension $(i+1)$.}

The sequence $C = (C^n)$ of chain complexes
is now a single chain complex in the abelian
category $\cm$ considered above (with an appropriately chosen growth rate $\mu$),
so we may apply the
construction of extended homology
and study the projective dimension, the torsion dimension, and the
Novikov - Shubin invariants. We will call these invariants {\it the asymptotic invariants}
of the sequence $K^n$. Note that the asymptotic invariants
really depend only on the geometry of $K^n$ for large $n\to \infty$. 

In order to construct the chain complex $C_\ast(K^n)$,
one has to choose orientations for all
simplices of $K^n$. But it is easy to see that different choices of orientations do not
influence the spectrum of the "Half - Laplacians" (4-13) and so the obtained
invariants do not depend on these orientations.

Note also that the asymptotic invariants are in general {\it geometric and not topological},
i.e. they will depend
on the simplicial decomposition of $K^n$'s and not on the topology of $K^n$.

For future references, let us make the following simple observation.

\proclaim{4.9. Proposition} Given a growth process $(K^n)$ as above, its asymptotic
invariants in dimension $i$ depend only on the growth process consisting of the skeletons
of $K^n$ of dimension $(i+1)$. In particular, the asymptotic invariants in dimension zero
depend only on the 1-skeletons of $K^n$.\qed
\endproclaim

\heading {\bf \S 5. Spectrum of towers: theorem of L\"uck}\endheading

In this section we will reformulate the theorem
of L\"uck \cite{L}.

L\"uck
considers a sequence of normal subgroups
$$\dots \subset \Gamma_{k+1}\subset \Gamma_k\subset\dots \subset\Gamma_1 \subset\pi$$
such that the index $[\pi:\Gamma_k]$ is finite for all $k$ and the intersection $\cap \Gamma_k$
is the trivial group. 
Let $X$ be a finite polyhedron with fundamental group
$\pi$. For each $k$ we have the finite sheeted covering $X^k\to X$ corresponding to the
subgroup $\Gamma_k$, and therefore we have a growth process $(X^k)$ (in the terminology
of section 4) determined by this tower of covering.
The theorem of L\"uck \cite{L} computes the asymptotic invariants of this growth process.

\proclaim{5.1. Theorem (Lu\"ck \cite{L})} Choose for the growth rate
$\mu = (\mu^k)$ the numbers
$$\mu^k = |\pi:\Gamma_k|^{-1}\tag5-1$$
(the inverses of the orders of the quotients $\pi/\Gamma_k$). Then

(i) The projective dimension of the growth process $(X^k)$
equals to the $L^2$ Betti number of the universal
covering of $X$ in the corresponding dimension.

(ii) The torsion dimension of the extended homology vanishes.
\endproclaim

We may conclude that the towers of coverings represent a very special class
of growth processes.

\heading {\bf  \S 6. Growing flat bundles}\endheading

Here we will consider an example of a growth process, which is a
generalization of the construction of tower of coverings, considered in the previous section.
We will fix a polyhedron $X$ and study a sequence of flat bundles over $X$ of growing
dimension. Our aim is to understand the asymptotic invariants in this situation.

This section contains only a general discussion of the problem;
the results are given by Theorems 8.2, 8.3, 8.4, 9.2 and 11.1.

\subheading{6.1} Let $X$ be a fixed finite simplicial polyhedron and let $\E^k$
be a sequence of finite dimensional flat bundles over $X$. We will assume that each bundle
$\E^k$ is supplied with a flat metric.

Define the growth rate $\mu = (\mu^k)$ as
$$\mu^k = (\rank \E^k)^{-1}.\tag6-1$$

For each integer $k$ we have the chain complex $C^k = C_\ast(X,\E^k)$ over $\C$.
The basis of this
chain complex is formed by the flat sections of $\E^k$ defined over the oriented simplices
of $X$. The boundary homomorphism is given by restricting a flat section over a simplex
$\sigma$
on all the faces of $\sigma$, multiplied by the sign, expressing compatibility of the
orientations of the simplex $\sigma$ with the orientation of the face.

We want to view this sequence of complexes $C^k$, where $k = 1,2,3,\dots$,
as a single complex in the category $\cm$. To meet all the requirements of section
4.6, we need to introduce
a scalar product in $C_\ast(X,\E^k)$. We will do it as follows: the scalar product of two
flat sections $s_1$ and $s_2$, which are defined
over two different simplices of $X$ is zero; if $s_1$ and $s_2$ are defined over
the same simplex $\sigma$ of $X$, then the scalar product $(s_1,s_2)$
equals to the scalar product of $s_1(v)$ and
$s_2(v)$ in the fiber of the bundle $\E^k$ over $v$, where $v$ is any vertex of the simplex
$\sigma$ - the result is independent on the choice of $v$, since the metric on $\E^k$
is supposed to be flat.

We obtain a chain complex $C = (C^k)$ in the abelian category $\cm$ and
we want to understand its asymptotic invariants.

\subheading{6.2} Note, that the construction of the tower of coverings (cf. \S 5)
is a special case of this construction.
In fact, in the situation of \S 5 for any $k$ we have the action of $\pi$ on the group
ring of the finite quotient $V^k = \C[\pi/\Gamma_k]$.
More precisely, we consider the action of $\pi$ from the left on the group algebra
$\C[\pi/\Gamma_k]$ and the corresponding flat bundle $\E^k$ over $X$.
Note that this bundle has a
flat metric, which comes from the metric of $\C[\pi/\Gamma_k]$ in which the elements of
$\pi/\Gamma_k$ form an orthonormal base.
The homology of the flat bundle $\E^k$ over $X$ coincides with the homology of the
normal covering $X^k\to X$, corresponding to $\Gamma_k$. 

\subheading{6.3. Example} Here we consider an example, which behaves unlikely the situation
with the towers of coverings.

Let $X$ be the closed 3-manifold obtained from the trefoil knot
\vskip 4cm
\centerline{Fig. 3}
by 0-framed surgery.
We have the canonical epimorphism $\phi: \pi_1(X) \to Z$
(the abelinization),
and therefore for any complex number $\xi$ with $|\xi| = 1$, there is a unique flat Hermitian
line bundle $\E_\xi$ with monodromy given by
$g \mapsto \chi_{\xi}(g) = \xi^{\phi(g)}$ for $g \in \pi$.
The dimension of homology $H_1(X, \E_\xi)$ is zero for all $\xi$ with $\xi^2 - \xi + 1 \ne 0$.
Here $\Delta(\xi) = \xi^2 - \xi + 1$ is the Alexander polynomial of the trefoil.
If $\xi$ is one of the
roots of the Alexander polynomial, i.e. if $\xi = \xi_\pm = e^{\pm \pi i/3}$,
then the dimension of the homology $H_1(X, \E_\xi)$ is 2.

Now, choose a sequence of complex numbers $\xi_k$ with $|\xi_k|=1$, such that $\xi_k\to \xi_+$.
Then we have a sequence of flat bundles $\E_{\xi_k}$, such that
{\it the sequence of dimensions $\dim H_1(X, \E_{\xi_k})$ may be an arbitrary sequence
consisting of
0 and 2}: we obtain 0 if $\xi_k \ne \xi_+$ and we obtain 2 if $\xi_k=\xi_+$.

Therefore, {\it
the projective dimension in this situation $\lo \dim_\C H_1(X, \E_\xi)$
may actually depend on the choice of the summation machine $\omega$}.

Note also that in this example
the corresponding characters $\chi_{\xi_k}$ converge and their limit is the
character $\chi_{\xi_+}$ at the root of the Alexander polynomial.

Suppose now that in the situation described above $\xi_k$ tends to $\xi_+$, but $\xi_k\ne \xi_+$.
Then we see that the projective dimension of the growth process is zero
(independently of $\omega$). However,
we will have the torsion dimension equal to 2. 

\heading{\bf \S 7. Characters of representations
and extended $L^2$-homology}\endheading

In this section we will show that the extended homology of a finite polyhedron $X$ with
coefficients in a representation $\M$ in a von Neumann category with a trace,
depends mainly on the
character $\chi_{\M}: \pi\to \C$
of the fundamental group of $X$ determined by $\M$. We also show that any
positive self-adjoint class function on the fundamental group can be realized as
the character of a unitary representation in a von Neumann category.

\subheading{7.1} Suppose that $\ca$ is a von Neumann category with a fixed trace $\tr$.

Let $\pi$ be a discrete group. We will consider representations of $\pi$ on objects of
$\ca$. More precisely, let $\M$ be an object of $\ca$; then a representation of $\pi$
is a ring homomorphism $\rho : \C[\pi] \to \hom_{\ca}(\M,\M)$. Such representation will
be called
{\it unitary} if $\rho$ is a $\ast$-homomorphism, i.e. if it preserves the involutions. Here we
assume that the group ring is supplied with the standard involution $g\mapsto g^{-1}$ for
$g\in \pi$. 

Any representation $\rho: \C[\pi] \to \hom_{\ca}(\M,\M)$ as above determines {\it the character} 
$$\chi_{\M}: \pi \to \C,\quad g\mapsto \tr(\rho(g)), \quad g\in \pi.\tag7-1$$

The character $\chi_{\M}$ is clearly constant on the conjugacy classes of $\pi$.
Also, if the representation is
unitary, then the character $\chi_{\M}$ has the property
$$\chi_{\M}(g^{-1}) = \overline{\chi_{\M}(g)}\tag7-2$$
for any $g\in \pi$. Class functions with
this property are called {\it self-adjoint}. Another important property of characters is
{\it positivity}: for any element $a\in \C[\pi]$ of the group algebra $\C[\pi]$ holds
$$\chi_{\M}(a^\ast a) \ge 0.\tag7-3$$

It is not true in general that the character
determines the representation up to the natural equivalence.

Using the construction of \cite{F}, we know that
to any finite CW space $X$ with fundamental group $\pi_1(X) = \pi$
we may assign extended homology $\H_\ast(X,\M)$ with coefficients in $\M$.

Our observation here is that (assuming that the trace $\tr$ on $\ca$ is normal)
the most important invariants of the extended homology
can be computed using only the character $\chi_{\M}$ of the representation $\M$:

\proclaim{7.2. Theorem} Suppose that the chosen trace $\tr$ on the von Neumann category
$\ca$ is normal. Let $X$ be a finite polyhedron with fundamental group $\pi$.
Then for any unitary representation $\rho: \pi \to \hom_{\ca}(\M,\M)$,
one can find the spectral density function $F_i(\lam)$ of the extended homology
$\H_i(X,\M)$ with coefficients in $\M$
using the character $\chi_{\M}$ of $\M$ as the only information on the representation $\M$.
In particular, the von Neumann dimension and the Novikov - Shubin invariants of
$\H_i(X,\M)$ depend only on the character $\chi_\M$ (and on $X$, of course).
\endproclaim

The proof of Theorem 7.2 is given in section 12.

We will see later in 10.1 that Theorem 7.2 is false without assuming normality of the
trace $\tr$.

Also, the Theorem is not true if the representation $\rho$ is not unitary.
Although Theorem 7.2 can be generalized to non-unitary representations, but the conclusion then
is different; we will consider this generalization elsewhere.

As a simple example, consider a finite dimensional unitary representation $V$ of $\pi$
and form the
tensor product $\M = V \otimes_\C \ell^2(\pi)$. The character $\chi_\M$ of this representation
equals to the character of $\M' = \ell^2(\pi)\oplus \dots \oplus \ell^2(\pi)$ ($\dim V$ times).
Then we obtain from Theorem 7.2 that the spectral density functions of $\H_i(X,\M)$ and
$\H_i(X, \M')$ coincide.

\subheading{7.3. Constructing representations with given characters} Here we will consider the
following problem: given a class function $\chi: \pi \to \C$, which is self-adjoint (7-2)
and non-negative (7-3), we want to construct a unitary representation
$\rho: \pi \to \hom_{\ca}(\M,\M)$ in certain von Neumann category $\ca$ with a normal trace
$\tr$ such that the character $\chi_{\M}$ of this representation is the given function $\chi$.
We will see that there is a canonical construction for this purpose. This construction
is very similar to the classical constructions (cf. \cite{N}, \S 30, and also \cite{G});
therefore we will be very brief.

First, we will associate a Hilbert space $\H_{\chi}$ with a given self-adjoint non-negative
class function $\chi: \pi \to \C$. We will denote by $J_{\chi}$ the following two-sided ideal
of $\C[\pi]$:
$$J_{\chi} = \{a\in \C[\pi]; \chi(ab) =0 \quad\text{for any}\quad b\in \C[\pi]\}.\tag7-4$$
Then we define the Hilbert space $\H_{\chi}$ as the completion of the factor-ring
$\C[\pi]/J_{\chi}$ with respect to the following scalar product
$$(a, b) = \chi(a b^\ast),\quad a,b \in \C[\pi].\tag7-5$$
It is easy to check that the obvious left and right actions of $\pi$ on the factor-ring
$\C[\pi]/J_{\chi}$ are continuous with respect to the norm determined by the scalar product
(7-5), and thus these actions extend to the left and right actions of $\pi$ on $\H_{\chi}$.
Both these actions are in fact unitary.

Note that the previous construction applied to the case when $\chi$ is the delta-function
at the unit element of the group $\pi$, gives the standard Hilbert space $\ell^2(\pi)$
which is usually associated with the group $\pi$.

Now we will construct a von Neumann algebra $\NN(\chi)$ acting on $\H_{\chi}$. We will denote
by $\NN({\chi})$ the space of all bounded linear maps $A: \H_{\chi}\to \H_{\chi}$, commuting
with the action of $\pi$ from the left. We obtain that
$$\C[\pi]/J_{\chi} \subset \NN(\chi)\tag7-6$$
(where $\C[\pi]$ acts from the right on the Hilbert space $\H_{\chi}$).

Now we will define the following function (the trace)
$$\tau: \NN(\chi) \to \C.\tag7-7$$
For $A\in \NN(\chi)$ set
$$\tau(A) = (A\cdot 1, 1)\tag7-8$$
where $1\in \C[\pi]/J_{\chi} \subset \NN(\chi)$ denotes the unit element and the brackets
$(\ ,\ )$ denote the scalar product (7-5). One easily check that:
\roster
\item $\tau$ {\it is a trace on the von Neumann algebra $\NN(\chi)$;
\item $\tau$ is normal;
\item $\tau$ is faithful;
\item on the subring $\C[\pi]/J_{\chi} \subset \NN(\chi)$
the trace $\tau$ coincides with $\chi$.}
\endroster

As shown in section 2.6 (example 3) of \cite{F}, the von Neumann algebra $\NN(\chi)$ acting on
$\H_\chi$ generates a finite von Neumann category $\ca$, where $\A$ is the group algebra
$\C[\pi]$. The trace $\tau$ on the algebra $\NN(\chi)$ determines a trace $\tr$ on the
category $\ca$. This trace on $\ca$ is clearly normal (since $\tau$ is normal).

Now, we have a unitary action of $\pi$ on $\M = \H_{\chi}\in\ob(\ca)$ and the
corresponding character $\chi_{\M}$ equals $\chi$.

\heading {\bf \S 8. Approximating characters}\endheading

Here we study the general problem about the relation between the von Neumann Betti numbers
and the dimensions of the homology of a sequence of finitely dimensional representations,
assuming that their characters converge to the character of the von Neumann representation.
We find a relation, involving an additional term, the torsion dimension,
which was studied in \S 2.

In the next section we consider the situation (which we call
arithmetic approximation) when this additional term vanishes.

\subheading{8.1} In this section we will study the following generalization of
the situation considered by W. L\"uck \cite{L}.

Suppose that $\pi$ is a discrete group and we are given a sequence of finite dimensional
unitary representations $\rho_k: \pi \to \End(V^k)$, where $k = 1, 2, \dots$. We will denote
by $\chi_k: \pi \to \C$ the corresponding characters. The dimensions
of these representations $\dim V^k = \chi_k(1)$ are not supposed to be constant. We will denote 
by $\mu^k = (\chi_k(1))^{-1}$ the inverse numbers; the numerical
sequence $\mu = (\mu^k)$ describes the growth rate
of the dimensions. We will consider also the {\it normalized characters}
$\tilde \chi_k = \mu^k\chi_k: \pi \to \C$, where $k = 1,2, \dots$.

Using the von Neumann category $\cm$ of section 4.1 (where
$\mu = (\mu^k)$ is the growth rate specified in the previous paragraph), we can view
the given sequence of representations $\rho_k: \pi \to \End(V^k)$, where $k=1,2,\dots$
as a single representation
$\rho_0: \pi \to \hom_{\cm}(\V,\V)$. Here $\V = (V^k)$ is the
object of $\cm$ determined by the given sequence $V^k$ of finite dimensional Hilbert spaces.

Using the construction of 1.8, for any finite polyhedron $X$ with fundamental group
$\pi$ we have the extended homology $\H_\ast(X,\V)$ with coefficients in $\V$. If we
choose the Dixmier type trace (4-4) on the category $\cm$, we obtain the numerical invariants
of the extended homology  -
the projective dimension, the torsion dimension, and the Novikov - Shubin invariants, cf. \S 1.
We will denote by $\dd P(\H_i(X,\V))$ the projective dimension and by
$\tdim_\omega T(\H_i(X,\V))$ the torsion dimension.
Recall that the projective dimension is just
$$\dd P(\H_i(X,\V)) = \lo [\frac{\dim_{\C} H_i(X, V^k)}{\dim_{\C} V^k}].\tag8-1$$
Note also that these invariants depend in general on the choice of the
summation procedure $\omega$. 

\proclaim{8.2. Theorem} In the situation described in subsection 8.1,
suppose that we are given another von Neumann category $\ca$ with a normal trace $\tr$
and a unitary representation $\rho: \pi \to \hom_{\ca}(\M,\M)$.
Suppose that the normalized characters of the finite dimensional
representations $\tilde \chi_k$ converge pointwise (as functions on the group $\pi$)
to the character $\chi_{\M}: \pi \to \C$ of $\M$, when $k\to \infty$. In other words,
we assume that for any $g\in \pi$ holds
$\lim_{k\to \infty} \tilde \chi(g) = \chi_{\M}(g).$
Then the following formula holds, which expresses the von Neumann dimension of
the extended homology $\H_i(X,\M)$
with coefficients in $\M$ (with respect to the trace $\tr$) by means of the dimensions of the
homology of the approximating finite dimensional representations:
$$\dim_{\tr}P(\H_i(X,\M)) = \dd P(\H_i(X,\V)) + \tdim_\omega T(\H_i(X,\V)).\tag8-2$$ 
\endproclaim
Note that the RHS of (8-2) contains only information obtained from the finite dimensional
flat bundles. The LHS of (8-2) is the $L^2$-Betti number with respect to a normal trace;
it is independent of the choice of the summation process $\omega$. Thus, (8-2) implies
that the sum of the projective dimension and the torsion dimension is independent of $\omega$.
Note, however, that the choice of $\omega$ may influence each of the numbers in the RHS of
(8-2), as example 6.3 shows. 

Compared with L\"uck's theorem \cite{L}, formula (8-2) contains an additional summand
(the torsion dimension). We will discuss in section 9 the conditions under which this
torsion dimension vanishes. 

The proof of Theorem 8.2 is given in \S 12.

The following statement shows that one may
recapture the entire spectral density function of the extended homology $\H_i(X,\M)$
in terms of the finite dimensional approximations.

\proclaim{8.3. Theorem (Approximation of the spectral density function)} Under the
condition of Theorem 8.2 the spectral density function $F_i(\lam)$ of extended homology
$\H_i(X,\M)$ and the spectral density function $G_i(\lam)$ of extended homology $\H_i(X,\V)$
are related as follows
$$F_i(\lam) = \lim_{\epsilon \to +0} G_i(\lam + \epsilon)\quad\text{for all}\quad
\lam \ge 0.\tag8-3$$
In other words, $F_i$ coincides with $G_i$, made right continuous.
\endproclaim

The proof of Theorem 8.3 is also postponed until section \S 12.

In the next theorem we point out conditions under which the projective
dimension and the torsion dimension in the RHS of (8-2) are both
independent of the summation procedure $\omega$.

\proclaim{8.4. Theorem (Analytic curve of representations)}
Let $\pi$ be a discrete group and let $\RR_N(\pi)$ denote the real
analytic variety of all representations of $\pi$ into the $N$-dimensi\-onal unitary group
$\U(N)$. Suppose, that
$\rho_n\in \RR_N(\pi)$, where $n = 1, 2, \dots$ is an infinite sequence of representations
such that there exists
a real analytic curve $\rho: [0,\epsilon) \to \RR_N(\pi)$ and a sequence
$t_n\in [0, \epsilon)$, such that $\rho_n = \rho(t_n)$, $t_n \to 0$, $t_n\ne 0$. Let
$\rho_0 = \rho(0)$ denote the limit representation. Then for any finite polyhedron $X$
and for any homomorphism $\phi : \pi_1(X) \to \pi$ we obtain the sequence of unitary
representations
$$\psi_n : \pi_1(X) @>\phi >> \pi @>{\rho_n}>> \U(N).\tag8-4$$
Each $\psi_n$ produces a flat $N$-dimensional unitary bundle over $X$, which we will
denote by $V^n$. Then the dimension of the homology $H_i(X, V^n)$
is constant for large $n$. Therefore, the torsion dimension of the sequence
$\V = (V^n)$ of these flat bundles equals to the jump in the Betti number
$$\tdim_\omega H_i(X,\V) = \dim H_i(X, V^0) - \dim H_i(X, V^n),\tag8-5$$
where $n$ is sufficiently large. In particular, we see that the torsion dimension
is independent of the choice of $\omega$. 
\endproclaim
\demo{Proof of Theorem 8.4} We only have to show that the dimension $\dim H_i(X, V^n)$
stabilizes for large $n$; the rest follows from Theorem 8.2.

We use the well known property of upper semi continuity of the dimension, cf.
\cite{H}, Ch.3, \S 12.
For any $t\in [0,\epsilon)$ denote by $V^t$ the flat bundle over $X$ with monodromy
$\phi\circ \rho(t)$. Then there exists a non-constant real analytic function $f(t)$
such that the dimension
$\dim H_i(X, V^t)$ assumes the constant value, say $D$, for all $t$ with $f(t) \ne 0$; moreover,
$\dim H_i(X, V^t) \ge D$ for all $t$. Suppose that we have a sequence of points $t_n$
with $t_n \to 0$,
$t_n \ne 0$. If $f(t_n)$ is zero only for finitely many $n$, then we obtain that the dimension
$\dim H_i(X, V^n)$ is constant for large $n$. However, if $f(t_n)=0$ for infinitely many $n$,
then the function $f(t)$ must be identically zero -- a contradiction. \qed
\enddemo

\heading{\bf \S 9. Arithmetic approximation}\endheading

It turns our that one may impose some arithmetic conditions on the approximating sequence
of flat bundles, which would imply vanishing of the additional correcting
term (the torsion dimension), appearing
in Theorems 8.2 and 8.4. Roughly, the arithmeticity condition requires that each
approximating finite
dimensional representation be definable over the ring of algebraic integers of an
algebraic number field, and the degrees
of these number fields must have a common upper bound.

This result implies the theorem of L\"uck.
Namely, L\"uck \cite{L} considers
the tower of finitely sheeted regular coverings, which is equivalent to studying
the homology of $X$ with coefficients in the representations
of $\pi$ on $\C[\pi/\Gamma_k]$ (cf. 6.2); these representations
are clearly defined over the integers $\Z$.

The main theorem of this section contains also a statement that
the torsion part of the extended homology $\H_i(X,\M)$
is of determinant class assuming that the character $\chi_\M$ of $\M$
admits an arithmetic approximation . This result generalizes a theorem
proven by D. Burghelea, L. Friedlander and
T. Kappeler \cite{BFK} for the case $\M = \ell^(\pi)$.
The proof presented in this paper (cf. \S 12),
is quite similar to the proof suggested in \cite{BFK}.

\subheading{9.1. Definition (Arithmetic approximation)} Suppose that $\chi:\pi \to \C$ is a
positive self-adjoint
class function $\chi: \pi \to \C$, cf. 7.1. We will assume that $\chi$ is normalized
so that
$\chi(1) =1$. We will say that $\chi$ {\it admits an arithmetic approximation},
if there exists a
sequence of finite dimensional unitary representations $\rho_k: \pi \to \End_{\C}(V^k)$,
where $k =1,2, \dots$, such that the following conditions are verified:

{\bf A}. Let $\tilde \chi_k: \pi \to \C$ denote the normalized character of $\rho_k$, i.e.
$\tilde \chi_k = \chi_k(1)^{-1}\chi_k$, where $\chi_k$ is the character of $\rho_k$.
Then the sequence $\tilde\chi_k(g)$ converges to
$\chi(g)$ for any $g\in \pi$.

{\bf B}. For each $k$ there is given an algebraic number field $\F_k\supset \Q$, imbedded into
the field of complex numbers $\C$ such that its image $\F_k\subset \C$ is invariant under
the complex conjugation. We will consider $\F_k$ together with the involution induced
from $\C$. We suppose that {\it the representation
$\rho_k$ can be defined over $\F_k$}. In other words, there exists a representation
$\tilde \rho_k : \pi \to
\End_{\F_k}(W_k, W_k)$, preserving a positively defined Hermitian form
$\<\ ,\ \>_k: W_k \times W_k \to \F_k\subset \C$,
which produces $\rho_k$ by extension of the scalars from $\F_k$ to $\C$.

{\bf C}. Denote by $\oo_k$ the ring of algebraic integers of $\F_k$. We suppose that for each
$k$ there exists an $\oo_k$-lattice $\L_k \subset W_k$
(i.e. a finitely generated $\oo_k$-submodule generating $W_k$ over $\F_k$),
which is invariant under the
action of $\pi$ and such that the form $\<\ ,\ \>_k$ restricted on $\L_k$ assumes values in
$\oo_k$.

{\bf D}. Denote by $\L_k^D$ the dual lattice
$$\L_k^D = \{w\in W_k; \<w,x\>\in \oo_k\quad\text{for any}\quad x\in \L_k\},$$
cf. \cite{FT}, page 122.
Then $\L_k\subset \L_k^D$ and the factor $ \L_k^D/  \L_k$ is a finite group.
We suppose that one can choose the lattices $\L_k$ such that
that there is an integer $M >0$ (independent of $k$) with
$M\cdot  \L_k^D/\L_k = 0$. 

{\bf E}. Let $h_k$ denote the degree of the number field $\F_k$ (cf. above) over the rationals.
We assume that there exists a common upper bound
$$h_k \le h\tag9-1$$
for all $k$.

{\bf F}. There exists a function $N: \pi \to \R_+$, having the property
$$N(gg')\le N(g)N(g')\quad\text{for any}\quad g,g'\in \pi,\tag9-2$$
and such that for any group element $g\in \pi$ the following
inequality holds
$$
|\sigma_j(\tilde \chi_k(g))| \le N(g)\tag9-3
$$
for all sufficiently large $k = 1, 2, ...$ and for all the embeddings
$\sigma_j: \F_k \to \C$, $j=1, \dots h_k$
of the algebraic number field $\F_k$.

Note, that conditions {\bf E} and {\bf F} imply that,
if the dimensions of the representations $\dim V^k$ are bounded, then
for any $g\in \pi$ the sequence $\tilde\chi_k(g)$ stabilizes for large $k$.
This statement
easily follows applying Lemma A (cf. subsection 12.3). However, this is not true if the
dimensions $\dim V^k$ grow; for example this stabilization does not happen in Theorem 0.5,
although all the conditions above hold.

The number $h$ (appearing in property {\bf E}) will be called {\it the degree} of the arithmetic
approximation. The number $M$ (appearing in {\bf D}) will be called {\it the denominator} of the
arithmetic approximation.

Now we show that the conditions of arithmetic approximation
imply vanishing of the torsion dimension, which appears in Theorem 8.2.

\proclaim{9.2. Theorem} Suppose that $\ca$ is a finite von Neumann category with a fixed
normal trace, $\M\in \ob(\ca)$ is a fixed object with $\dim_{\tr}(\M) =1$,
and $\rho: \pi \to \End_{\ca}(\M)$ is a unitary representation, having character
$\chi_{\M}: \pi \to \C$.
Suppose that we are given a sequence of finite dimensional representations $\rho_k:\pi \to
\End_\C(V^k)$ which provide an arithmetic approximation of the character $\chi_{\M}$, cf.
section 9.1.  Then for any finite polyhedron $X$ with fundamental group $\pi$ holds:

(i) the sequence of the normalized Betti numbers
$$\frac{\dim_{\C} H_i(X,V^k)}{\dim_\C V^k}$$
converges and its limit coincides with the von Neumann dimension
of the projective part of the extended homology
$$\lim_{k\to\infty} \frac{\dim_{\C} H_i(X,V^k)}{\dim_\C V^k} = \dim_{\tr} \H_i(X,\M).\tag9-4$$

(ii) Let $F_i(\lam)$ denote the spectral density function of extended homology
$\H_i(X,\M)$. Then the following inequality holds
$$F_i(\lam) - F_i(0) \le \frac{c}{-\log(\lam)}\tag9-5$$
for small $\lam >0$, where $c>0$ is a constant. 
(iii) For any $i$ the torsion part of the extended homology $\H_i(X,\M)$ is of
determinant class.
\endproclaim

L\"uck's theorem \cite{L} follows from this by taking $V^k=\C[\pi/\Gamma_k]$,
cf. 6.2. This flat bundle
can be defined over the integers. 

Intuitively, the integrality condition in Theorem 9.2 allows to conclude at some point
of the proof, that certain nonzero quantity cannot be too small.

The proof of Theorem 9.2 is given in section 12.

For the definition of the notion {\it determinant class} (which appears in the statement (iii)
of Theorem 9.2)
we refer to \cite{BFKM}. Cf. also \cite{CFM},
section 3.8, where it is explained why the condition of being of determinant class depends
only on the torsion part of the extended homology.

It is natural to ask for which groups $\pi$ the character $\chi$
of the natural representation of
$\pi$ on $\ell^2(\pi)$ (which is the delta-function at the unit $1\in \pi$) admits an
arithmetic approximation. It is easy to see that {\it it happens if and only if $\pi$
is residually finite.}
In fact, if $\pi$ is residually finite, we may construct the arithmetic
approximation as in L\"uck's theorem: if $\pi\supset \Gamma_1\supset \Gamma_2\dots $
is a  sequence of normal subgroups with trivial intersection, then we may take
$V^k = \C[\pi/\Gamma_k]$, which can be realized over integers. Conversely, if we are given
an arithmetic approximation (cf. 9.1)
with $\tilde \chi_k$ converging to $\chi$, then for any $g\in\pi$,
where $g\ne 1$, we have $\tilde \chi_k(g)\to 0$, and so there exists $k$ with
$\tilde\chi_k(g)<1$. This implies that the image of $g$ under the $k$-th representations
$\rho_k$ is nontrivial, $\rho_k(g)\ne 1$. Since the automorphism of the lattice
$\GL_{\oo_k}(\L_k)$ is residually finite, we obtain that $\pi$ must be residually finite.

Now we will give proofs of Theorems 0.3, 0.4 and 0.5 (cf. \S 0), deduced from Theorem 9.2.

\subheading{9.3. Proof of Theorem 0.3} 
We will use the notations introduced
in Theorem 0.3.

With the sequence of subgroups $\pi \supset \Gamma_1\supset \Gamma_2 \dots$
we associate the sequence of the
unitary representations $\rho_k: \pi \to \End(\Z[\pi/\Gamma_k])$ defined over $Z$. Here
we assume that $\pi$ acts on the group ring $\Z[\pi/\Gamma_k]$ as the left
regular representation. It is clear that if $V^k$ denotes the flat bundle over $X$ determined
by this representation then $H_i(X,V^k)\simeq H_i(\tilde X_k)$.

All conditions of arithmetic approximation (of \S 6.1)
are obviously satisfied. We only need to compute the normalized character $\tilde \chi_k$
of $\rho_k$. An elementary calculation shows that for $g\in \pi$
$$\tilde \chi_k(g) = \frac{n_k(g)}{n_k},$$
where $n_k$ is the total number of different subgroups of $\pi$ conjugate to $\Gamma_k$,
and $n_k(g)$ is the number of them, containing $g$.
Therefore, our assumption (0-1) implies that the normalized characters $\tilde \chi_k$
converge pointwise to the character of the standard representation of $\pi$ on $\ell^2(\pi)$.
Applying Theorem 9.2, we complete the proof. \qed

\subheading{9.4. Proof of Theorem 0.4}
Consider the representation
$$
\nu_k: \pi \to \End_{\oo_k}(\Z[\pi]\otimes_{\Z[\Gamma_k]} \oo_k^{m_k}),\tag9-6
$$
induced by the given unitary representation
$$
\rho_k: \pi \to \End_{\oo_k}(\oo_k^{m_k}).\tag9-7
$$
We want to apply Theorem 9.2 to the obtained 
sequence of representations $\nu_k$. It is clear that if $X$ is any polyhedron with
$\pi_1(X)=\pi$, then the homology of $X$ with coefficients twisted by $\nu_k$ is the same as the
homology of the covering space $\tilde X_k$ with coefficients twisted by $\rho_k$
(here $\tilde X_k \to X$ denotes covering corresponding to $\Gamma_k$).

>From the general properties of induced representations (cf. \cite{CR}, \S 10) we
see that the conditions {\bf B, C, D, E}
of section 9.1 are satisfied. We only need to check that the sequence of
normalized characters of $\nu_k$
converge to $\chi_0 : \pi \to \C$, where $\chi_0(g) = 0$ for $g\in \pi, g\ne 1$ and
$\chi_o(1)=1$. Note that $\chi_0$ is the character of the standard representation of
$\pi$ on $\ell^2(\pi)$ with respect to the von Neumann trace. Also, we need to check condition
{\bf F} in 9.1.

If we denote by $\chi_k: \Gamma_k \to \C$ the character of $\rho_k$ and by
$\eta_k: \pi \to \C$ the character of $\nu_k$; then
$$\eta_k(g) = \cases \chi_k(g)\quad \text{for}\quad g\in \Gamma_k,\\
0,\quad\text{if}\quad g\notin \Gamma_k.\endcases
$$
Therefore we see that for any $g\ne 1$ the normalized character $\tilde \eta_k(g)$ vanishes
for all large $k$. Thus {\bf F} of section 9.1 holds with $N(g)=0$.
Applying Theorem 9.2 finishes the proof. \qed

\subheading{9.5. Proof of Theorem 0.5} Let $\chi: \pi \to \oo_{\F}$ denote the character
of the given representation $\rho: \pi \to \Mat(m\times m, \oo_{\F})$. Then the character
of the tensor power $\rho^{\otimes k}$ is $g \mapsto \chi(g)^k$. We claim that for $g \ne 1$
holds $|\chi(g)| < m$ and therefore the normalized character of $\rho^{\otimes k}$
tends to 0:
$$\frac{|\chi(g)|^k}{m^k} \to 0.$$
In fact, $\rho(g)$ viewed as a complex $m\times m$ matrix, can be diagonalized, and on the
diagonal we will obtain $m$ numbers with norm 1. Therefore $|\chi(g)|\le m$ and the equality
holds if and only if $\rho(g)$ belongs to the center.

Thus we obtain condition {\bf A} of 9.1. Conditions {\bf B, C, D, E} are obvious.
Condition {\bf F} follows from the assumption that the field $\F$ is cyclotomic: then all
the Galois transformations preserve the complex norm.
\qed

Next we will formulate a corollary of Theorem 9.2 which may be useful.

\proclaim{9.6. Theorem} Let $\F\subset \C$ be an algebraic number field invariant under the
complex conjugation, and let $\F' \subset \C$ be a cyclotomic field. We will denote by
$\oo_{\F}$ and $\oo_{\F'}$ the corresponding rings of algebraic integers. Let
$\pi$ be a discrete group and let $\rho: \pi \to \Mat(n\times n,\oo_{\F})$ and
$\rho_k: \pi \to \Mat(n_k\times n_k,\oo_{\F'})$, where $k = 1, 2, \dots$,
be unitary representations,
such that for any $g\in \pi$ the limit
$$\lim_{k\to \infty} \frac{\chi_k(g)}{n_k} = \chi_0(g)\tag9-8$$
exists; here $\chi_k: \pi \to \C$ denotes the character of $\rho_k$. Let $X$ be a compact
polyhedron with $\pi_1(X) = \pi$ and let $\E$ and $\E_k$ (for $k=1, 2, \dots$) denote the
complex flat vector bundles over $X$ determined by $\rho$ and $\rho_k$ correspondingly.
Then the sequence of the normalized Betti numbers
$$\frac{\dim H_i(X, \E\otimes \E_k)}{n\cdot n_k}, \quad\text{where}\quad k=1, 2, \dots\tag9-9$$
converges and its limit can be found as follows.
Let $\ca$ be a finite von Neumann category with a normal trace $\tr$ and let $\M$ be an
object of $\ca$ supplied with a unitary action of $\pi$ having the character
$$g \mapsto \chi_{\M}(g) \ = \ \frac{\chi_0(g)\chi(g)}{n},\quad g\in \pi,\tag9-10$$
where $\chi$ denotes the character of $\rho$ (we know from \S 7 that such von Neumann
representation $\M$ exists). Then the limit of the sequence (9-9) equals to the dimension
(with respect to the trace $\tr$) of the extended homology $\H_i(X, \M)$:
$$\lim_{k\to \infty} \ \frac{\dim H_i(X, \E\otimes \E_k)}{\dim \E\cdot \dim \E_k} \ = \
\dim_{\tr} \H_i(X, \M).\tag9-11$$
\endproclaim

\demo{Proof}Theorem 9.6 follows by applying Theorem 9.2 similarly
to the arguments given in 9.3, 9.4, 9.5. We will only point out here
how one constructs the function
$N: \pi \to \R_+$, which appears in {\bf F}. For $g\in \pi$ we define
$$N(g) = n^{-1}\cdot \sup_j\{||\sigma_j(\chi(g))||\},\tag9-12$$
where $\sigma_j: \F \to \C$ runs over all the embeddings of $\F$.
We consider the representations
$\rho\otimes\rho_k$, where $k = 1, 2, \dots$,
as defined over the ring of algebraic integers of the
compositum $\F\F'$ of the fields $\F$ and $\F'$; any embedding of
$\F\F'$ into $\C$ determines embeddings of $\F$ and $\F'$, and (9-12) is clearly
enough to establish property {\bf F} (cf. 6.1), since $\F'$ is assumed to be cyclotomic
and so its Galois transformations are unitary.
\qed
\enddemo

\heading{\bf \S 10. Examples}\endheading

\subheading{10.1} Here we show that Theorem 7.2 is false if the trace $\tr$ on von
Neumann category $\ca$ is not normal.

We will construct two finite von Neumann categories with traces (one normal and one not normal)
and two unitary representations of the fundamental group of a polyhedron $X$ on objects of these
categories, such that the characters of this representations are equal but the projective
dimensions of the corresponding homology are distinct.

As the first von Neumann category we will take the category $\Cal C_1$ of finite dimensional
Euclidean vector spaces with the usual trace.
As the second category $\Cal C_2$ we will take the category
$\cm$ (cf. \S 4), where $\mu$ is the constant sequence $\mu^k =1$. We will consider the Dixmier
type trace $\tr_\omega$ in $\Cal C_2$, cf. (4-4).

Now we will return to the
example described in subsection 6.3. The space $X$ was obtained by 0-framed surgery on
the trefoil knot, and for any $\xi\in S^1$ we had a unitary flat bundle $\E_\xi$ over $X$.
We will suppose that the sequence of 
points $\xi_k$ on the unit circle is chosen so that $\xi_k\to \xi_+$ and $\xi_k\ne \xi_+$,
where $\xi_+=e^{\pi i/3}$ is a root of the Alexander polynomial, cf. 6.3.
Then the sequence of flat bundles $\E_{\xi_k}$, viewed as a single flat bundle $\V$ with fiber
an object of $\Cal C_2$, has character $g \mapsto \xi_+^{\phi(g)}$, where $g\in \pi_1(X)$
and $\phi: \pi_1(X) \to \Z$ is the abelinization.
We see that the same character has the line bundle $\E_{\xi_+}$ (viewed as bundle in $\Cal C_1$).
Then we have $\dim_{\tr_{\omega}} \H_1(X,\V) = 0$, however $\dim H_1(X,\E_{\xi_+}) = 2$. 

\subheading{10.2} Here we will show that Theorem 9.2 is false without assumption (9-1) that
the degrees of the number fields $h_k$ are bounded.

We will again use the example 6.3. We choose  points $\xi_k$ on the unit circle such that
$\xi_k$ converges to $\xi_+$ and $\xi_k$ for any $k$ is a root of 1. We will assume that
$J = \{k; \xi_k =\xi_+\}$ is a subsequence; we may actually choose the sequence $\xi_k$
such that $J$ is arbitrary. Note that the corresponding sequence of flat line bundles
$\E_{\xi_k}$ satisfies all the conditions of arithmetic approximation, cf. 9.1, besides (9-1).
We see that the sequence of dimensions $\dim H_1(X,\E_{\xi_k})$ is the following: we have 2,
for $k\in J$ and we have 0, for $k\notin J$. Thus sequence (9-4) is not convergent.

\subheading{10.3. Algebraic integers on the unit circle} Here (preparing tools for the next
example) we observe that {\it there
exist algebraic integers on the unit circle, which are not roots of unity.}
The simplest example is as follows. Consider the roots of the equation
$$
\aligned
&z^4 - z^3 - z^2 - z +1 \ = \\
= \ &(z^2 - \frac{1+\sqrt{13}}{2}\cdot z + 1)\cdot
(z^2 - \frac{1-\sqrt{13}}{2}\cdot z + 1) \ = \ 0.
\endaligned\tag10-1
$$
Two of its roots are complex, lying on the unit circle, they are roots of the second factor in
(10-5).
We will denote them $e^{i\alpha}$ and $e^{-i\alpha}.$ Here $\alpha \simeq 130.6463$ degrees.

Two other roots are real, we will denote them by $r$ and $r^{-1}$, where $r \simeq 0.5807$.
\vskip 4cm
\centerline{Figure 4.}
The numbers $e^{i\alpha}$, $e^{-i\alpha}$, $r$ and $r^{-1}$ are algebraic integers, which are all
conjugate to each other.
We conclude that $e^{i\alpha}$ is not root of 1 since otherwise all its conjugates
would be roots of 1, and so they would be points of the unit circle.

Also, the numbers $e^{i\alpha}$, $e^{-i\alpha}$, $r$ and $r^{-1}$
are in fact units of the corresponding ring of algebraic integers.

We observe that the powers
of  $e^{i\alpha}$ are dense on the unit circle. The powers of $r$ tend to 0.

\subheading{10.4} Here we will show that Theorem 9.2 is false without condition {\bf F} in 9.1.

Let $X$ be the 3-manifold obtained by 0-framed surgery from the trefoil, as in example 6.3.
We will use the notations introduced in 6.3 and in 10.1 and 10.2.

Let $e^{i\alpha}$ denote the algebraic integer on the unit circle, constructed in 10.3. The
powers $e^{in\alpha}$, where $n \in \Z$, are dense on the circle, and thus we can find a
subsequence $n_k$ such that $e^{in_k\alpha}$ converges to $\xi_+$ (recall that $\xi_+$ denotes
the root of the Alexander polynomial of the trefoil). We will denote by $\E_k$ the unitary flat
line bundle over $X$ corresponding to the point $e^{in_k\alpha}$ (as in 10.1). Then we obtain,
that the sequence of Betti numbers
$\dim H_1(X,\E_k)$ consists of zeros, and the corresponding characters
converge to the character of $\E_{\xi_+}$, but $\dim H_1(X,\E_{\xi_+})=2$.

Note that in this example all
the conditions of arithmetic approximation of 9.1 except {\bf F} are satisfied.
Our field $\F$ in this example has 4 embeddings $\sigma_j: \F \to \C$, $j=1,2,3,4.$
The embedding, which sends the number $e^{i\alpha}$ to $r^{-1}$ (cf. notations of 10.3)
sends $e^{in_k\alpha}$ to $r^{-n_k}$, which tends to $\infty$, violating {\bf F}.

\heading{\bf \S 11. Approximation in characteristic $p$}\endheading

\proclaim{11.1. Theorem} Let $p$ be a prime number. Suppose that
$$\pi \supset \Gamma_1 \supset \Gamma_2 \supset \dots, \quad\text{where}\quad
 \cap \Gamma_j = \{1\},\tag11-1$$
is a chain of normal subgroups such that for each $j$ the index $[\pi : \Gamma_j]$ is a power
of $p$. Let $X$ be a finite $CW$ complex with fundamental group $\pi$ and let $\tilde X_j \to X$
be the normal covering corresponding to $\Gamma_j$. Then for any $i$ the sequence
$$\frac{\dim_{\fp}H_i(\tilde X_j,\fp)}{[\pi : \Gamma_j]},\quad j=1,2, \dots\tag11-2$$
decreases and so the limit
$$\lim_{j\to \infty}\frac{\dim_{\fp}H_i(\tilde X_j,\fp)}{[\pi : \Gamma_j]}\tag11-3$$
exists.
\endproclaim
\proclaim{11.2. Corollary} If for some prime $p$ the fundamental group $\pi$ of a finite
CW complex $X$ admits a
chain of normal subgroups (11-1) such that all the factors $\pi/\Gamma_j$
are $p$-groups, then the following
inequality holds
$$b_i^{(2)}(X) \le \dim_{\fp} H_i(X,\fp),\tag11-4$$
for the $L^2$-Betti number $b_i^{(2)}(X)$.
\endproclaim
Corollary 11.2 follows immediately from Theorem 11.1 using the Theorem of L\"uck
\cite{L} and the inequality
$$\dim_{\C} H_i(\tilde X_j, \C) \le \dim_{\fp}H_i(\tilde X_j,\fp).$$
\demo{Proof of Theorem 11.1} Using Corollary on page 25 of \cite{La}, we may assume that
in the given chain
of normal subgroups (11-1) all the factors $\pi/\Gamma_j$ are cyclic of order $p$.
Thus, to prove Theorem 11.1 it is enough to show that
$$\dim_{\fp} H_i(\tilde X_{j+1}, \fp) \le p\cdot \dim_{\fp}H_i(\tilde X_j,\fp)\tag11-5$$
for any $j$.

Fix a tringulation of $X$ and consider the induced triangulations on all the coverings
$\tilde X_{j}$.

>From this moment we will assume that $j$ is fixed. We will consider the $p$-sheeted covering
$\tilde X_{j+1} \to \tilde X_j$.
Let $C$ denote the chain complex of simplicial chains of $\tilde X_{j+1}$ with coefficients
in the finite field $\fp$. $C$ is a free
finitely generated chain complex over the ring $\Lambda = \fp[\Z/p]$. Note that $\Lambda$
has a unique maximal ideal $\m = (t-1)\Lambda$, where $t$ denotes
the generator of $\Z/p$. We have
the following filtration on $\Lambda$:
$$\Lambda \supset \m \supset \m^2\supset \dots \supset \m^{p-1} \supset 0.$$
Therefore, we obtain the filtration
$$C \supset \m C \supset \m^2 C\supset \dots \supset \m^{p-1}C \supset 0$$
and all the factor-complexes $\m^r C/\m^{r+1}C$, where $r=0, 1, \dots p-1$,
are isomorphic to the chain complex of
$\tilde X_j$ with coefficients in $\fp$. We obtain that there is a spectral sequence,
starting from
$$\bigoplus_{p \ \text{times}} H_i(\tilde X_j,\fp)$$
and converging to $H_i(\tilde X_{j+1},\fp)$. This proves (11-5). \qed

\enddemo

\subheading{11.3. Questions}

Can limit (11-3) be greater than the $L^2$-Betti number $b^{(2)}_i(X)$?

Does sequence (11-3) always stabilize after a finitely many steps?

\heading{\bf \S 12. Proofs of Theorems 7.2, 8.2, 8.3, 9.2}\endheading

Here we finally present proofs of the main theorems of this paper.
These proofs are related to each other. Therefore
we use the same notations and terminology. In fact, we assume that the reader will read
the proofs in the proper order (7.2, 8.2, 8.3 and then 9.2 - lexicographical ordering!).
Also, we very much use arguments of L\"uck's paper \cite{L}, and sometimes we do not repeat them,
but instead refer to \cite{L}.  Thus, it will be very helpful for the
reader to have a copy of \cite{L} at hand while reading this section. 

\subheading{12.1. Proof of Theorem 7.2}
Suppose that $X$ has a fixed tringulation. Consider the
chain complex $C_\ast(\tilde X)$ of the simplicial chains in the universal covering $\tilde X$.
It is a
complex of free finitely generated $\Z[\pi]$-modules; its basis is formed by the lifts of the
oriented simplices of $X$. Note that each chain module $C_i(\tilde X)$ is naturally
supplied with a non-degenerate
$\Z[\pi]$-valued Hermitian scalar product which is defined using the basis formed by the
lifts of the cells as the orthonormal basis. The boundary homomorphism
$d: C_{i+1}(\tilde X) \to C_i(\tilde X)$ is given by the matrix with entries in $\Z[\pi]$.
Consider the "adjoint" homomorphism $d^\ast: C_i(\tilde X) \to C_{i+1}(\tilde X)$ which is
defined using the above mentioned $\Z[\pi]$-valued Hermitian scalar product. Then we have
the following self-adjoint homomorphism
$$d^\ast d : C_{i+1}(\tilde X) \to C_{i+1}(\tilde X), \quad
d^\ast d \in \Z[\pi] \otimes \Mat(a\times a, \Z),\tag12-1$$
(where $a$ denotes the number of $(i+1)$-dimensional simplices in $X$).
If $p(z)$ is any polynomial with real coefficients then we may form $p(d^\ast d)$ and the
result will be a self-adjoint matrix with entries in $\R[\pi]$.
Now, applying the character $\chi_{\M}$ to this matrix produces a matrix $\chi_{\M}(p(d^\ast d))$
with entries in $\C$, which is Hermitian. We will consider then the trace (in the usual sense)
of this Hermitian matrix $\Tr(\chi_{\M}(p(d^\ast d)))$. Note that the same answer will be
obtained if we will first map the matrix $d^\ast d$ via the representation
$\rho:\C[\pi] \to \Hom_{\ca}(\M,\M)$,
then applying the polynomial
$p(z)$ to get
$$p(\rho(d^\ast d)) \in \Hom_{\ca}(\M^{a}, \M^{a}),\tag12-2$$
and finally computing the trace
$\tr_{\M^a}$ of (12-2):
$$\Tr(\chi_{\M}(p(d^\ast d))) = \tr_{\M^{a}} (p(\rho(d^\ast d))).\tag12-3$$
This follows from the definition of the trace on a category (cf. 1.3) and the definition of the
character.

We would like to be able to compute (using only the character)
more general expressions of the form
$\Tr(\chi_{\M}(f(d^\ast d)))$, where $f(z)$ is
a real valued function. The most important for us is the case, when the function
$f(z)$ above is the
characteristic function of an interval $[0,\lam^2]$, which we will denote by $f_\lam(z)$.

According to W. L\"uck \cite{L}, this can be done as follows.
Choose a sequence of real polynomials $p_n(z)$
such that
$$p_n(z) \to f_\lam(z)\quad\text{and}\quad |p_n(z)|\le L,\tag12-4$$
where both properties (12-4) hold for any $z\in [0,N]$. Here $N$ is a fixed apriori large 
number such that
$$|\rho(d^\ast d)|\le N\quad\text{for any unitary representation}\quad \rho.\tag12-5$$
We will take $N$ to be $a$ times the sum of all coefficients, which appear in the
matrix elements of $d^\ast d$ (this claim is similar to Lemma 2.5 in \cite{L}). To be more
precise, we know, that $d^\ast d = (b_{ij})$, where the entries $b_{ij}$ of this
$a\times a$-matrix belong to the group ring $\C[\pi]$, $b_{ij}= \sum \beta_{ij}(g)\cdot g$,
where the sum is taken over $g\in \pi$ (only finitely many terms are nonzero). We define $N$ as
$$N = \sum_{i,j,g} |\beta_{ij}(g)|.$$

Using the Lebesgue theorem on Majorized convergence and the assumption that the trace
$\tr$ is normal and, therefore it is continuous with respect to the ultraweak topology on
$\hom_{\ca}(\M^{a}, \M^{a})$ (cf. \cite{Di}, Part I, chapter 6, \S 1), we obtain that
the operator
$p_n(\rho(d^\ast d))$ converges ultraweakly to
$$f_\lam(\rho(d^\ast d)) = \int_0^{\lam^2}\lam dE_\lam,\tag12-6$$
in the von Neumann algebra $\hom_{\ca}(\M^{a}, \M^{a})$,
where $E_\lam$ is the right continuous
spectral decomposition of $\rho(d^\ast d)$. Thus using (12-3),
we find
$$\Tr(\chi_{\M}(p_n(d^\ast d))) =
\tr_{\M^{a}}(p_n(\rho(d^\ast d))) \to \tr_{\M^{a}}(f_\lam(\rho(d^\ast d))).\tag12-7$$
We obtain finally the following formula for the spectral density function $F_i(\lam)$ of the
extended homology $\H_i(X,\M)$:
$$F_i(\lam) = \lim_{n\to \infty} \Tr(\chi_{\M}(p_n(d^\ast d))),\quad \lam > 0.\tag12-8$$
The last formula involves only the character $\chi_{\M}$.
Since $F_i(\lam)$ is right continuous, we
find also (using only the character $\chi_{\M}$)
the von Neumann dimension of the extended homology
$\dim_{\tr} P(\H_i(X, \M))$ as the limit
$\lim_{\lam\to +0} F_i(\lam)$. This completes the proof of Theorem 7.2.
\qed

\subheading{12.2. Proof of Theorems 8.2 and 8.3}
The proof uses the methods of L\"uck \cite{L} with certain adjustments.
We will use the notations introduced in the proof of Theorem 7.2, cf. 12.1.
In particular we will use
formula (12-8). As in the proof of 7.2 we will denote by $F_i(\lam)$
the spectral density function of the extended homology $\H_i(X,\M)$. Since the trace $\tr$
is assumed to be normal, we will assume that $F_i(\lam)$ is right continuous. The von Neumann
dimension $\dim_{\tr}P(\H_i(X,\M)$ is by the definition $F_i(0)$. 

For $k=1,2, \dots$ denote by $F^k_i(\lam^2)$ the spectral density function of the finite
dimensional operator
$\rho_k(d^\ast d)$, where $\rho_k: \pi \to \End(V^k)$ is the $k$-th representation.
As in the proof of 7.2, we regard here $d^\ast d$ as the matrix with entries in the group
ring $\Z[\pi]$, i.e. $d^\ast d \in \Z[\pi] \otimes \Mat(a\times a,\Z)$, where $a$
denotes the number of $(i+1)$-dimensional cells in $X$. Therefore,
$\rho_k(d^\ast d) \in \End(V^k)\otimes \Mat(a\times a,\Z)$. Similarly
to (12-8) we have for $\lam >0$
$$F_i^k(\lam) = \lim_{n\to \infty}\Tr( \chi_k(p_n(d^\ast d))),\tag12-9$$
where $p_n(z)$ is any sequence of polynomials constructed as in the proof of Theorem 7.2.
Let us introduce also the functions
$$
\align
&G_i(\lam) = \lo [\mu^k F^k_i(\lam)]\quad\text{and}\tag12-10\\
&G^+_i(\lam) = \lim_{\epsilon\to +0} G_i(\lam +\epsilon),\tag12-11
\endalign
$$
defined for $\lam\ge 0$.
According to our definitions, we have
$$
\align
&G_i(0) = \dd P(\H_i(X,\V))\qquad\text{and}\tag12-12\\
&G^+_i(0) =
\dd P(\H_i(X,\V)) + \tdim_\omega T(\H_i(X,\V)).\tag12-13
\endalign
$$
Therefore to prove Theorem 8.2 we have to show that $F_i(0) = G_i^+(0).$

Note that $G_i(\lam)$ is the spectral density function of extended homology
$\H_i(X,\V)$ as defined in \cite{F}, section 3.7.
We will see now that the second function $G_i^+(\lam)$ is
in fact more important.

We will choose the polynomials $p_n(z)$ as follows.
Denote by $g_n: \R \to \R$ the function
$$
g_n(z) = \cases
1+ 1/n\quad\text{for}\quad z\le \lam^2, \\
1 + 1/n - n(z-\lam^2)\quad\text{for}\quad \lam^2 \le z \le \lam^2 +1/n,\\
1/n \quad\text{for}\quad \lam^2 +1/n \le z
\endcases\tag12-14
$$
and construct the polynomials $p_n(z)$ such that
$$g_n (z) \le p_n(z) \le 2 \quad\text{and}\quad \lim_{n\to \infty}p_n(z) = f_\lam(z)$$
for all $z\in [0,N]$. Here $f_\lam(z)$ denotes the characteristic function of the interval
$[0,\lam^2]$ and $N$ is the large number constructed in the proof of Theorem 7.2, cf. (12-5).

With this choice of the polynomials $p_n(z)$ we may show
that for any $n$, $k$ and $\lam >0$ holds
$$\mu^k F^k_i(\lam) \le \Tr(\tilde\chi_k (p_n(d^\ast d))) \le (1+1/n)\mu^k F_i^k(\lam + 1/n)
 + a/n,\tag12-15$$
where $a$ denotes the number of $(i+1)$-dimensional cells in $X$. To prove this one denotes
by $E_k(\lam)$ the ordered set of eigenvalues $z$ of $\rho_k(d^\ast d)$ satisfying $z\le \lam$
listed with multiplicities. Then
$$\Tr(\tilde \chi_k(p_n(d^\ast d))) = \mu^k\cdot\sum_{z\in E_k(\lam^2)} p_n(z),$$
and now to obtain (12-15) one just repeats the arguments on page 469 of \cite{L}. 

Taking in (12-15) for fixed $n$ the limit $\lo$ with respect to $k$ and using the assumption
that $\tilde \chi_k \to \chi_{\M}$ we obtain
$$G_i(\lam) \le \Tr(\chi_{\M}(p_n(d^\ast d))) \le (1+1/n) G_i(\lam +1/n) + a/n.\tag12-16$$
Therefore, taking the limit in (12-16) when $n\to \infty$ and using (12-8) we get
$$G_i(\lam) \le F_i(\lam) \le G^+_i(\lam).\tag12-17$$
>From the last inequality we obtain for $\epsilon >0$
$$F_i(\lam) \le G^+_i(\lam) \le G_i(\lam +\epsilon) \le F_i(\lam +\epsilon),\tag12-18$$
and since we know that $F_i(\lam)$ is right continuous,
this shows (by passing to the limit when $\epsilon \to 0$) that
$$F_i(\lam) = G_i^+(\lam).\tag12-19$$
This is precisely the statement of Theorem 8.3.

Since both functions in the last equality are right continuous, we obtain $F_i(0) = G_i^+(0)$,
which completes the proof of Theorem 8.2 (cf. 12-13)). \qed

\subheading{12.3. Proof of Theorem 9.2.(i) and 9.2.(ii)} We will use the following Lemma
from algebraic number theory:

\proclaim{Lemma A} Let $\F_k\subset \C$ be a number field of degree $h_k \le h$ and
let $\oo_k$
be the ring of algebraic integers of $\F_k$. Let
$\sigma_1, \sigma_2, \dots, \sigma_{h_k}: \F_k \to \C$
denote all the distinct embeddings of $\F_k$ into the complex numbers. Then for any element
$\alpha\in \oo_k$ with $\alpha\ne 0$, the condition
$$|\sigma_i(\alpha)| \le R\quad\text{for all}\quad i=1, 2, \dots, h_k\tag12-20$$
implies
$$|\sigma_j(\alpha)| \ge R^{1-h_k}\tag12-21$$
for any $j = 1, 2, \dots, h_k$.
\endproclaim

This Lemma is well known, however we will give a simple independent proof.
Similar argument is used in \cite{Sh}, in the proof of Theorem 11 in chapter 1.

\demo{Proof of Lemma A} The product
$$\prod_{i=1}^{h_k}\sigma_i(\alpha)$$
(the norm of $\alpha$) is a nonzero integer. Therefore we obtain
$$|\sigma_j(\alpha)| \ge \prod_{i=1, i\ne j}^{h_k} |\sigma_i(\alpha)|^{-1} \ge R^{1-h_k}.$$
This completes the proof of Lemma A.
\enddemo

Here is another lemma, which we will need:
\proclaim{Lemma B} Let $A=(a_{ij})$ be a $k\times k$-matrix with complex entries. Suppose that
for some $C>0$ and $K\ge 1$ holds
$$|\Tr(A^r)| \le C\cdot K^r\tag12-22$$
for all $r = 1, 2, \dots k$. Then we have the following estimate for the coefficients
$s_r = s_r(A)$ of the characteristic polynomial
$\det(\lam - A) = \sum_{r=0}^k (-1)^{k-r}s_{k-r}(A)\lam^r$
of $A$:
$$|s_r(A)| \le \frac{C(C+1)\dots (C+r-1)}{r!}\cdot K^r.\tag12-23$$
\endproclaim
\demo{Proof of Lemma B} For $r=1, 2, \dots k$ denote
$$p_r = \Tr(A^r)= \sum_{i=1}^k \lam_i^r,$$ 
where $\lam_1, \lam_2, \dots \lam_k$ denote the eigenvalues of $A$. We have 
$$s_r = s_r(A) = \sum_{i_1<\dots < i_r} \lam_{i_1}\lam_{i_2}\dots \lam_{i_r}.\tag12-24$$

We will prove (12-23)
by induction on $r$ using the following Newton's identity
$$(-1)^r\cdot s_r = s_1p_{r-1} - s_2p_{r-2} + \dots +(-1)^r s_{r-1}p_1 - p_r,\tag12-25$$
cf. \cite{CR}, page 314.

Since $s_1 = p_1$, the inequality (12-23) holds for $r=1$; suppose that it has been
established for all values of $r$ which are less than the given one. Then from (12-25) we obtain
$$r\cdot |s_r| \le C\cdot K^r \cdot \{\sum_{j=1}^{r-1} \binom {C+j-1}j +1\},$$
and now the desired inequality (12-23) follows from the identity
$$\frac{C}{r}\cdot \{\sum_{j=1}^{r-1} \binom {C+j-1}j +1\} \ = \ \binom {C+r-1}r,\tag12-26$$
which can be easily checked by induction. This completes the proof of Lemma B.
\enddemo

Now we will prove statement (i) of Theorem 9.2.
We will use the notations introduced in the proofs of
Theorems 7.2 and 8.2. Also we will use the notations introduced in 9.1.

Let us fix some $k$. We know that the matrix
$\rho_k(d^\ast d) \in \End(V^k)\otimes \Mat(a\times a, \Z)$ is congruent over $\C$ to the matrix
$\tilde \rho_k(d^\ast d) \in \End_{\F_k}(W^k)\otimes \Mat(a\times a, \Z)$
(by condition {\bf B} in 9.1. The latter matrix has entries in the field $\F_k$
and therefore the characteristic polynomial
$$q_k(t) = \det(t - \rho_k(d^\ast d))\tag12-27$$
(of $\rho_k(d^\ast d)$ or equivalently of $\tilde\rho_k(d^\ast d)$)
has coefficients in $\F_k$. Write
$q_k(t) = t^{\nu}\overline q_k(t)$, where $\overline q_k(0)\ne 0$ and
$\overline q_k(0)\in \F_k$.

We claim now that
$$M^{a\dim V^k}\cdot \overline q_k(0)\quad\text{belongs to}\quad \oo_k,\tag12-28$$
where $M$ is the
denominator of the arithmetic approximation, cf. 9.1. To show this we note, that from
our assumptions {\bf C} and {\bf D} in 9.1 it follows that $M$ times the dual (over the field
$\F_k$) of the
$\oo_k$-homomorphism
$$d: \L_k\otimes _\pi C_{i+1}(\tilde X) \to \L_k\otimes _\pi C_{i}(\tilde X)$$
is well defined as a homomorphism
$$ \L_k\otimes _\pi C_{i}(\tilde X) \to \L_k\otimes _\pi C_{i+1}(\tilde X).$$
In other words,
$$M\cdot \tilde \rho_k(d^\ast d): W_k \otimes_\pi C_{i+1}(\tilde X) \to
W_k \otimes_\pi C_{i+1}(\tilde X)$$
preserves the $\oo_k$-lattice $\L_k \otimes_\pi C_{i+1}(\tilde X)$. If we assume
that the last lattice
is free over $\oo_k$, then we may take its basis as the basis of
$W_k\otimes_\pi C_{i+1}(\tilde X)$ and compute the characteristic polynomial with respect
to this basis. We obtain that the polynomial
$M^{a\dim V^k}\cdot q_k(tM^{-1})$ has all coefficients in $\oo_k$. Therefore
all coefficients of $M^{a\dim V^k-\nu}\cdot \overline q_k(tM^{-1})$ belong to $\oo_k$
and therefore $M^{a\dim V^k-\nu}\cdot \overline q_k(0)\in \oo_k$, which implies our
statement (12-28).

In order to prove (12-28) in the general case (not assuming that the lattice
$\L_k \otimes_\pi C_{i+1}(\tilde X)$ is free), we proceed as follows. For any prime ideal
$\p\subset \oo_k$ we consider the localization of $\L_k \otimes_\pi C_{i+1}(\tilde X)$
with respect to the complement of $\p$ which is now a free $(\oo_k)_\p$-module
(since $(\oo_k)_\p$ is a principal ideal ring, cf. \cite{FT}, page 59). The
arguments of the previous paragraph show that the valuation
$$v_\p(M^{a\dim V^k}\cdot \overline q_k(0))\ge 0$$
is non-negative. Since this is true for any prime ideal $\p$ of $\oo_k$, we obtain
that the fractional ideal generated by $M^{a\dim V^k}\cdot \overline q_k(0)$
is contained in $\oo_k$, which proves (12-28).

We know that
$\overline q_k(0)$ is the product of all the nonzero eigenvalues of the matrix
$\rho_k(d^\ast d)$
and from (12-5) we know that any of these eigenvalues is less or equal than $N$,
where the number
$N\ge 1$ was constructed
in the proof of Theorem 7.2. Note that $N$ is determined only by the matrix $d^\ast d$
(i.e. by the polyhedron $X$)
and does not depend on $k$. Therefore we obtain that
$$|\overline q_k(0))| \le N^{a\dim V^k}, \quad i=1,2,\dots, h_k.\tag12-29$$

We claim that a similar estimate 
$$|\sigma_j(\overline q_k(0))| \le N_1^{a\dim V^k}, \quad i=1,2,\dots, h_k. \tag12-30$$
holds for any embedding $\sigma_j: \F_k \to\C$
of the number field $\F_k$ into $\C$, where the constant $N_1$ is independent on $k$ and of
$j = 1, 2, \dots h_k$.
To show this, we will denote by $g_1, g_2, \dots g_s$ all elements of $\pi$, which appear
with nonzero coefficients in the matrix $d^\ast d$. We set
$$N_1 = 4 \cdot N \cdot L,\quad\text{where}\quad L = \max\{N(g_1), N(g_2),\dots, N(g_s), 1\};
\tag12-31$$
the numbers $N(g_i)$ are given by property {\bf F} in 9.1. In order to prove inequality (12-30)
we will apply Lemma B to the following $(a\dim V^k \times a\dim V^k)$-matrix
$$A= \sigma_j(\rho_k(d^\ast d))\tag12-32$$
for fixed $k$ and $j$; here we view $\rho_k(d^\ast d)$ as
the matrix with entries in the field $\F_k$ and we obtain a complex matrix applying the
embedding $\sigma_j: \F_k \to \C$. An obvious argument using our assumption (9-3) gives
the estimate
$$|\tr(A)| \le a\dim V^k \cdot N\cdot L.\tag12-33$$
Similarly, using our assumption (9-4) about the behavior of the function $N(g)$, gives
$$|\tr(A^i)| \le a\dim V^k \cdot (N L)^i\tag12-34$$
for all $i$. Now, by Lemma B, we may conclude (using the notation introduced in Lemma B) that
$$|s_r(A)| \le \binom {a\dim V^k +r -1}{r}\cdot (NL)^r\tag12-35$$
for $r= 0, 1, 2, \dots, a\dim V^k$. Clearly, $\sigma_j(\overline q_k(0) = s_r(A)$ for
the largest $r\le a\dim V^k$ with $s_r(A)\ne 0$, and since we have the following obvious
estimate for the binomial coefficient
$$\binom {a\dim V^k +r -1}{r} \ \le \ 2^{2a\dim V^k -1} \ < \ 4^{a\dim V^k},\tag12-36$$
combining (12-35) and (12-36), we obtain (12-30).

Using inequality (12-21) of Lemma above and also (12-23), we obtain
$$|M^{a\dim V^k}\overline q_k(0)| \ge (MN_1)^{a\dim V^k (1-h_k)},\tag12-37$$
which implies
$$|\overline q_k(0)| \ge M^{-ha\dim V^k} N_1^{(1-h)a\dim V^k}.\tag12-38$$

Now we use Lemma 2.8 of L\"uck  \cite{L}. We apply it to the operator $\rho_k(d^\ast d)$
and estimate (12-38).
We obtain (using the
notations introduced in the proofs of Theorems 7.2 and 8.2 and after some elementary
transformations)
$$\mu^k[F^k_i(\lam) - F^k_i(0)] \le \frac{c}{-\ln(\lam)},\tag12-39$$
where the constant $c$ is
$c = h\log M + (h-1)\log N_1 + \log N.$
Taking $\lo$ with respect to $k\to \infty$ in (12-27) we obtain
$$G_i(\lam) - G_i(0) \le \frac{c}{-\ln(\lam)}.\tag12-40$$
By the definition
$\tdim_\omega T(\H_i(X,\V)) = G^+_i(0) - G_i(0)$
(cf. (12-12) and (12-13)). Comparing this with
(12-40) we see that the torsion dimension vanishes.

Now to finish the proof we only have to show that the sequence $\mu^k\dim H_i(X,V^k) =
\mu^k F^k_i(0)$ converges. To do so, we will introduce the following notations:
$$\underline G_i(\lam) = \lim\inf_{k\to \infty} \mu^k F^k_i(\lam),$$
$$\overline G_i(\lam) = \lim\sup_{k\to \infty} \mu^k F^k_i(\lam),$$ 
$$\underline G_i^+(\lam) = \lim_{\epsilon \to +0} \underline G_i(\lam +\epsilon),$$
$$\overline G_i^+(\lam) = \lim_{\epsilon \to +0} \overline G_i(\lam +\epsilon).$$
>From (12-15) we obtain (by passing to the limits with respect to $k$)
$$\overline G_i(\lam) \le \Tr(\chi_{\M}(p_n(d^\ast d)))) \le
(1+1/n)\underline G_i(\lam + 1/n) + a/n$$
and the limit $n\to \infty$ gives
$$\overline G_i(\lam) \le F_i(\lam) \le \underline G_i^+(\lam).$$
Now we see that for $\epsilon >0$ we have
$$F_i(\lam) \le \underline G_i^+(\lam) \le \underline G_i^+(\lam +\epsilon) \le
\overline G_i(\lam +\epsilon) \le F_i(\lam +\epsilon)$$
and thus when $\epsilon \to 0$ we get
$$F_i(\lam) = \underline G_i^+(\lam) = \overline G_i^+(\lam).\tag12-41$$
In particular, we obtain
$$F_i(0) = \underline G_i^+(0) = \overline G_i^+(0).\tag12-42$$

On the other hand, using inequality (12-27) we obtain
$$\underline G_i(\lam) \le \underline G_i(0) + \frac{c}{-\ln(\lam)},$$
$$\overline G_i(\lam) \le \overline G_i(0) + \frac{c}{-\ln(\lam)},$$
which give for $\lam \to 0$
$$\underline G_i^+(0) = \underline G_i(0),\qquad\text{and}\qquad
\overline G_i^+(0) = \overline G_i(0).$$
Comparing the last equalities with (12-42) gives $\underline G_i(0) = \overline G_i(0)$.
This proves the convergence of $\mu^k\dim H_i(X,V^k)$ and completes the proof of (i).

Statement (ii) was proven above by (12-40).

\subheading{12.4. Proof of Theorem 9.2.(iii)}
The arguments here are similar to those used by D.Burghelea, L.Friedlander and
T.Kappeler, cf. Appendix of \cite{BFK}.

We will use the notations introduced in the proof of Theorem 8.2 and in
Theorem 7.2. We want to show that
$$ \ld (\rho(d^\ast d)) \ \overset {def}\to = \
\int_{+0}^\infty \ln(\lam)dF_i(\lam) >-\infty.\tag12-43$$
If $N$ is the large number constructed in Proof of Theorem 7.2, (cf. (12-5)),
then $F_i(\lam)$ is
constant for $\lam \ge N$ and so we may write (integrating by parts)
$$
\aligned
\ld (\rho(d^\ast d)) =&  (F_i(N) - F_i(0))\cdot \ln N \ +\\
+&\lim_{\epsilon \to +0}\{[F_i(\epsilon) - F_i(0)](-\ln \epsilon) -
\int_\epsilon^N \frac{F_i(\lam)- F_i(0)}{\lam}d\lam\}.
\endaligned\tag12-44
$$
>From the last formula we obtain the inequality
$$\ld (\rho(d^\ast d)) \ge  (F_i(N) - F_i(0))\cdot \ln N \ - 
\int_{+0}^N \frac{F_i(\lam)- F_i(0)}{\lam}d\lam.\tag12-45
$$
Note that in the similar formula for the finite dimensional operator $\rho_k(d^\ast d)$ we
have the equality
$$\ld (\rho_k(d^\ast d)) =  (F_i^k(N) - F_i^k(0))\cdot \ln N \ - 
\int_0^N \frac{F_i^k(\lam)- F_i^k(0)}{\lam}d\lam\},\tag12-46
$$
since the spectral density function $F_i^k(\lam)$ is constant for small $\lam \ne 0$.
Now using inequality (12-38) we find
$$
\ld (\rho_k(d^\ast d)) = \log |\overline{q}_k(0)| \ge -ha\log(MN_1)\cdot \dim V^k,\tag12-47
$$
where $a$ denotes the number of $(i+1)$-dimensional simplices in $X$ and $N_1$ is given by
(12-31). Multiplying by
$\mu^k = (\dim V^k)^{-1}$, we get
$$\mu^k\cdot \ld (\rho_k(d^\ast d)) \ge  -ha\log(MN_1) .\tag12-48$$
Similarly to statement 1 of Lemma 3.3 of L\"uck \cite{L}, we obtain the inequality
$$\int_0^N \frac{F_i(\lam)- F_i(0)}{\lam}d\lam \le \lim \inf_{k\to \infty}
\{ \int_0^N \frac{\mu^k F_i^k(\lam)- \mu^k F_i^k(0)}{\lam}d\lam \}.\tag12-49$$
Now, we multiply equality (12-34) by $\mu^k$ and pass to the limit infimum, when $k$ tends to
infinity. Since $\mu^kF_i(N) \to F_i(N)$ (by Theorem 8.3) and $\mu^kF_i(0) \to F_i(0)$
(by Theorem 9.2.(i)), we obtain finally (combining inequalities (12-45), (12-48) and (12-49))
$$\ld (\rho(d^\ast d)) \ge - ha\log (MN_1).\tag12-50$$
This completes the proof. \qed

\enddocument